\documentclass[11pt]{article}

\usepackage{amsmath}
\usepackage{amssymb}

\newcommand{\Prob}{\operatorname{P}}

\newcommand{\R}{{\mathbf R}}

\newcommand{\Lloc}{L^1_{\mathrm{loc}}}
\newcommand{\Mloc}{M_{\mathrm{loc}}}
\providecommand{\cite}[1]{[#1]}

\newcommand{\sca}[2]{\langle #1, #2\rangle}

\DeclareMathOperator{\divv}{div}

\DeclareMathOperator{\supp}{supp}

\makeatletter
\newcommand{\footn}%
{
\renewcommand{\@makefnmark}{} 
\footnotetext{\mbox{\tiny\jobname 
    \quad \the\year\ \the\month\ \the\day} \hfill \hfill}
}

\let\0 = \emptyset

\let\p=\partial
\let\e=\varepsilon
\let\lb=\lambda

\let\de=\delta
\let\G=\Gamma
\let\g=\gamma
\let\vf=\varphi 
\let\rar=\rightarrow
\let\a=\alpha 
\let\b=\beta

\let\Om=\Omega
\let\ti=\widetilde
\let\sm=\setminus
\let\ss=\subset
\let\s=\sigma

\let\w=\omega

\let\i=\infty

\let\hat=\widehat 
 
\let\bar=\overline

\title{ Support theorems for the Radon transform and \\ Cram\'er-Wold theorems}
\author{Jan Boman
  \footnote{Department of Mathematics, Stockholm University, 
    SE-10691 Stockholm,
    Sweden. E-mail: jabo@math.su.se}, 
  Filip Lindskog
  \footnote{Department of Mathematics, Royal Institute of Technology, 
    SE-10044 Stockholm, 
    Sweden. E-mail: lindskog@kth.se}}

\numberwithin{equation}{section}                                    

\begin{document}

\maketitle

\begin{abstract}
This article presents extensions of the Cram{\'e}r-Wold theorem to 
measures that may have infinite mass near the origin.
Corresponding results for sequences of measures are presented 
together with examples showing that the assumptions imposed are sharp.
The extensions build on a number of results and methods concerned
with injectivity properties of the Radon transform.
Using a few tools from distribution theory and Fourier analysis
we show that the presented injectivity results for the Radon transform lead to 
Cram{\'e}r-Wold type results for measures. 
One purpose of this article is to contribute to making known to
probabilists interesting results for the Radon transform that have 
been developed essentially during the 1980ies and 1990ies.
\end{abstract}




\section{Introduction}

  The Cram\'er-Wold theorem \cite[p.~291]{CW36} states that a probability 
  measure $P$ on $\R^d$ is uniquely determined by the values it gives 
  to halfspaces $H_{\omega,p}=\{x\in\R^d;\, x\cdot\omega<p\}$ for 
  $\omega \in S^{d-1}$ and $p\in\R$. Equivalently, $P$ is uniquely 
  determined by its one-dimensional projections 
  $P\pi_{\omega}^{-1}$, where $\pi_{\omega}$ is the projection
  $\R^d \ni x \mapsto x\cdot\omega \in \R$ for $\omega \in S^{d-1}$.
  Moreover, a sequence of probability measures $P_k$ converges weakly
  to a probability measure $P$ in the sense that
  $\lim_{k\to\infty}\int\varphi\, dP_k=\int\varphi\, dP$ for all bounded 
  continuous real-valued $\varphi$, if for every $\omega \in S^{d-1}$
  $\lim_{k\to\infty}P_k(H_{\omega,p})=P(H_{\omega,p})$ for all but at most 
  countably many $p \in \R$. 
In recent years there has been an interest in analogous theorems in situations where the measure $P$ is not necessarily a probability measure but may have infinite mass near the origin. 
Such measures arise for instance as limits of scalings
of probability measures in multivariate extreme value theory 
(see e.g.~\cite{R87}) and limit theorems for sums of random 
vectors (see e.g.~\cite{MS01});
other examples include L\'evy measures
for infinitely divisible probability distributions and intensity 
measures for random measures (see e.g.~\cite{DV88} and \cite{S99}).
If $P$ has infinite mass near the origin, the value of $P(H_{\w,p})$  is of course not defined when the closure of $H_{\w,p}$ contains the origin, and the problem therefore becomes to decide if $P$ is determined by its values on all closed halfspaces contained in $\R^d\sm\{0\}$. 
We present three types of extensions of the Cram\'er-Wold theorem in that direction, and we show by examples that the assumptions made cannot be omitted. 
Our measures may take both positive and negative values unless
the contrary is explicitly stated.

For a finite (signed) measure $\mu$ with a density 
$f \in L^1(\R^d)$  we can write 
\begin{align}\label{halfspacerepr}
  \mu(H_{\omega,p})=\int_{-\infty}^p\left(\int_{L_{\omega,r}}f ds\right)dr,
\end{align}
where $ds$ is the Euclidean surface measure on the hyperplane 
$L_{\omega,p}=\{x\in\R^d;\,x\cdot\omega=p\}$.
The inner integral in \eqref{halfspacerepr}
is identified as the Radon transform of $f$ evaluated 
at the hyperplane $L_{\omega,r}$. Cram\'er-Wold theorems are therefore equivalent to injectivity  theorems for the Radon transform.  

The extension of the theory of the Radon transform to {\it distributions} $f$ is easy and well known \cite{H80}. Measures are distributions of order zero, and we can therefore without difficulty form the Radon transform of the measures we need to consider. 
In particular, an analogue of (\ref{halfspacerepr}) for measures is (\ref{dpmuh}) below. 
Since distributions are defined as linear forms on spaces of test functions, it is natural for us to define measures as linear forms on a space of continuous test functions (see Section 2).  
Working with measures in the way we do requires a very small part of distribution theory. In a few cases we shall use distributions of higher order than zero. A few facts from distribution theory that may not be well known to our readers are collected in an appendix.

In Section 2 we introduce measures as distributions of order zero and define the Radon transform and other operations on measures. 

In Section 3 we present four injectivity theorems for the Radon transform, here called Theorems A -- D. Theorems B -- D treat the case --- often called the exterior Radon transform --- when a function (measure) is to be reconstructed outside a compact, convex set $K$, the case which is in focus in this article. Theorem A, the injectivity theorem for the standard Radon transform, is included for completeness. Theorem B is the well known Helgason support theorem for the Radon transform; here uniqueness is guaranteed by the assumption that the measure is {\it rapidly decaying} at infinity (see definition below). That the rapid decay assumption cannot be omitted is well known: for any integer $m\ge d$ there exist functions $f$, homogeneous of degree $-m$, with Radon transform $Rf(L)=0$ for all hyperplanes $L$ not containing the origin. On the other hand, if $f$ is homogeneous of non-integral degree, then $f$ is uniquely determined by its exterior Radon transform (Theorem C). Theorem D, finally, proves injectivity for the exterior Radon transform if the unknown measure is supported in a closed, convex cone containing no complete straight line. 

In Section 4 we begin by presenting four Cram\'er-Wold type uniqueness theorems, Theorems 1 -- 4, parallel to Theorems A -- D, respectively. Theorem 3 occurs in two variants, Theorem 3a and Theorem 3b; in the latter case the unknown measure is assumed to be non-negative  and therefore does not have to be assumed homogeneous. The main part of Section 4 presents four Cram\'er-Wold theorems for sequences of measures, Theorems 1$'$ -- 4$'$. The main novelties of our paper are probably Theorem 4 and Theorem 4$'$, which are analoguous to Theorem D. 
To some extent Theorems 3b and 3$'$ are perhaps also new, although a similar 
result has been shown by Basrak, Davis and Mikosch \cite{BDM02}. 
The support assumption in Theorems 4  and 4$'$ are often satisfied in applications with the cone $Q$ being the positive orthant $\{(x_1,\ldots,x_d)\in\R^d;\, x_k\ge 0\}$. A particular case concerns sequences of scalings of probability measures (Corollary 2); 
this result answers  affirmatively the conjecture in \cite{BDM02}.

Examples showing that the assumptions in Theorem B -- D are sharp are given in Section~5. A shorter description of essentially the same examples appeared in \cite{Bo92}, page 28. Those examples show immediately that the assumptions in Theorems 2 -- 4 and 2$'$ -- 4$'$ are sharp. Moreover, choosing $f(x) = q(|x|)h(x)$, where $h$ is a non-trivial solution to $Rh(\w,p) = 0$ in $p\ne 0$ and $q(|x|)$ is a very slowly oscillating radial function, shows that the assertion of Corollary 2 is not true if $\b$ is an integer. A similar example was previously given by Hult and Lindskog \cite{HL061}.

As explained above,
this article presents extensions of the Cram{\'e}r-Wold theorem to 
measures that may have infinite mass near the origin.
The extensions build on a number of results and methods, that have 
been developed essentially during the 1980ies and 1990ies, concerned
with injectivity properties of the Radon transform.
One purpose of this article is to contribute to making known to 
probabilists interesting results for the Radon transform that 
have appeared in the mathematical literature.
As is well known, the Radon transform and its generalizations have 
been studied extensively after the invention of Computerized Tomography
in the 1970ies. 
Using a few tools from distribution theory and Fourier analysis
we show that the presented injectivity results for the Radon transform lead to 
Cram{\'e}r-Wold type results for measures. 
The paper is self-contained and only a minimum of basic tools from 
distribution theory and Fourier analysis are needed.
In particular, we have aimed to convince the reader of the usefulness
of some very basic facts from distribution theory for treating problems 
occurring in applications of probability theory.

\section{Measures and their Radon transforms}
Let $C_0(\R^d)$ be the space of continuous functions on $\R^d$ that tend to zero at infinity, equipped with the supremum norm $\|\cdot\|$. The dual space of  $C_0(\R^d)$, 
the space of continuous linear forms on $C_0(\R^d)$, will be denoted $M(\R^d)$. This is the space of signed measures with finite total mass. 
Throughout the paper a measure is
real-valued, as opposed to non-negative, unless anything else is said.
The action of a linear form $\mu$ on the test function $\vf$ will be denoted $\sca{\mu}{\vf}$, and the
norm of $\mu \in M(\R^d)$ as a linear form will be denoted $\|\mu\|_M$, that is, 
$$
\|\mu\|_M = \sup\{|\sca{\mu}{\vf}|;\, \vf \in C_0(\R^d), \ \|\vf\| \le 1 \} . 
$$
This is the total mass, or total variation norm $|\mu|(\R^d)$, of the measure $\mu$. The space $L^1(\R^d)$ is identified with a subspace of $M(\R^d)$ by $f\in L^1(\R^d)$ being identified with the linear form $\vf \mapsto \int f(x) \vf(x) dx$. 
The relationship between a measure
$\mu \in M(\R^d)$ considered as a set function $\mu(E)$ defined 
on the family of Borel
sets and $\mu$ considered as a linear form 
$\vf \mapsto \langle\mu,\vf\rangle$
is well known and explained by
the Riesz Representation Theorem 
(see e.g.~Theorem 6.19 in \cite{R66}), 
which says that to every continuous linear form $\Phi$ on $C_0(\R^d)$
there corresponds a unique Borel measure $\mu$ such that
\begin{align}\label{rrt}
  \langle\Phi,\vf\rangle=\int_{\R^d}\vf\, d\mu,
  \quad \vf\in C_0(\R^d),
\end{align}  
and $\|\Phi\|_M=|\mu|(\R^d)$.
The integral in \eqref{rrt} is defined in any 
textbook on measure and integration theory.
Conversely, if the set function $\mu(E)$ is given, then it is clear that \eqref{rrt}
defines a continuous linear form on $C_0(\R^d)$.

Any element $\mu\in M(\R^d)$ can be uniquely extended as a linear form to the space $C_b(\R^d)$ of bounded continuous functions on $\R^d$. This is perhaps most easily seen using the expression $\int \vf\, d\mu$ considering $\mu$ as a set function. If $\mu$ is considered as a linear form on 
$C_0(\R^d)$, we take a compactly supported continuous function $\chi$ that is equal to $1$  in some neighborhood of the origin and define 
\begin{equation}    
\sca{\mu}{\vf} = \lim_{A\rar\i} \sca{\mu}{\chi(\cdot/A)\vf(\cdot)} , \quad \vf \in C_b(\R^d) . 
\end{equation}
It is easy to see that this definition is independent of the choice of $\chi$.

Following Laurent Schwartz \cite{S66} we shall denote by $\mathcal D(\R^d)$ the space of $C^{\i}$ functions with compact support. Note that $\mathcal D(\R^d)$ is dense in $C_0(\R^d)$. 

For $f\in L^1(\R^d)$ the Radon transform $Rf$ is defined by 
$Rf(L) = \int_L f\, ds$, where $ds$ is the Euclidean surface measure on the hyperplane $L$, or 
\begin{equation}      \label{Rf}
Rf(\w,p) = \int_{L_{\w,p}} f \, ds , \quad (\w,p) \in S^{d-1}\times\R ,
\end{equation}    
where ${L_{\w,p}}$ is the hyperplane $\{x\in\R^d;\, x\cdot\w = p\}$. Note that $Rf$ is even,
$Rf(\w,p) = Rf(-\w,-p)$, since $L_{\w,p} = L_{-\w,-p}$. 
If $f\in L^1(\R^d)$, then $Rf$ is defined almost everywhere on $S^{d-1}\times\R$, and in fact $Rf(\w,\cdot)$ is in  $L^1(\R)$ for every $\w\in S^{d-1}$.  It is clear that 
$\|Rf(\w,\cdot)\|_{L^1(\R)} \le \|f\|_{L^1(\R^d)}$ for every $\w$ and that 
$\|Rf\|_{L^1(S^{d-1}\times\R)} \le \|f\|_{L^1(\R^d)}$. Here the norm in $L^1(S^{d-1}\times\R)$ is defined using the normalized surface measure on $S^{d-1}$, which we denote by $d\w$.

For $\mu\in M(\R^d)$ we shall define the Radon transform $R\mu$ as a  measure on $S^{d-1}\times\R$, that is, as a linear form on the space $C_0(S^{d-1}\times\R)$ of continuous functions on $S^{d-1}\times\R$ that tend to zero at infinity. To do this we need some more notation. If $\phi\in L^1(S^{d-1}\times\R)$ and $\psi\in C_0(S^{d-1}\times\R)$ we write 
$$
\sca{\phi}{\psi} = \int_{S^{d-1}} \int_{\R} \phi(\w,p) \psi(\w,p)  dp\, d\w .
$$
More generally, if $\phi$ is a linear form on $C_0(S^{d-1}\times\R)$ we write 
$\sca{\phi}{\psi}$ to denote the action of $\phi$ on the test function $\psi$; 
thus $\phi\in L^1(S^{d-1}\times\R)$ is identified with the linear form
$C_0(S^{d-1}\times\R)\ni\psi\mapsto\sca{\phi}{\psi}$. 
The dual Radon transform $R^*$ is defined for $\psi \in C_0(S^{d-1}\times\R)$ by 
\begin{equation}      \label{R*}
R^*\psi(x) = \int_{S^{d-1}} \psi(\w, x\cdot\w) d\w  .
\end{equation}     
If $\psi$ is even,  $\psi(\w,p) = \psi(-\w,-p)$, then $\psi$ can be considered as a function on the manifold of hyperplanes, $\psi(L_{\w,p}) = \psi(\w,p)$, and the geometric meaning of (\ref{R*}) is that $R^*\psi(x)$ is defined as the mean of $\psi(L)$ taken over all hyperplanes $L$ 
containing $x$.  It is easy to verify that $R^*$ is the adjoint of $R$ in the sense that 
\begin{equation}   \label{adjoint}
\sca{R\phi}{\psi} = \sca{\phi}{R^*\psi}
\end{equation}
for $\phi\in L^1(\R^d)$ and $\psi\in C_0(S^{d-1}\times\R)$. 
Therefore it is natural to define the 
Radon transform $R\mu$ of $\mu\in M(\R^d)$ by
\begin{equation}      \label{Rmu1}
\sca{R\mu}{\psi} = \sca{\mu}{R^*\psi} ,\quad \psi\in C_0(S^{d-1}\times\R) .
\end{equation}     
It is obvious that $R^*$ maps $C_0(S^{d-1}\times\R)$ and 
$C_b(S^{d-1}\times\R)$
into $C_b(\R^d)$, but in fact 
 $R^*$ maps $C_0(S^{d-1}\times\R)$ into $C_0(\R^d)$:
 
\smallskip
\noindent
{\bf Lemma 1.}
If $\psi \in C_0(S^{d-1}\times\R)$ then  $R^*\psi\in C_0(\R^d)$ and
\begin{equation}      \label{supR*}
\sup|R^*\psi| \le \sup|\psi| . 
\end{equation}   

\noindent
{\it Proof.} 
It is obvious that $R^*\psi$ is continuous and bounded and that (\ref{supR*}) holds, so we only need to prove that   $R^*\psi(x) \rar 0$ as $|x|\rar \i$. First observe that for any $A$ the measure of the set 
$$
E(x,A) = \{\w\in S^{d-1};\, |x\cdot\w| < A\} 
  =  \{\w\in S^{d-1};\, |(x/|x|)\cdot \w| <  A/{|x|}\}
$$
tends to zero as $|x|\rar \i$. Choose $A$ so that $|\psi(\w,p)| < \e$ for $|p| > A$, and then choose $B$ so large that the measure of $E(x,A)$ is less than $\e$ if $|x| > B$. Then, if $|x|>B$,  
\begin{equation*}
|R^*\psi(x)|  \le \big|\int_{E(x,A)}  \psi(\w, x\cdot\w) d\w \big| 
      + \big|\int_{\complement E(x,A)}  \psi(\w, x\cdot\w) d\w \big|   
       \le  \e  \sup|\psi|  + \e ,
\end{equation*}
which completes the proof.

It follows from the definition that $R\mu \in M(S^{d-1}\times\R)$, the space of measures on $S^{d-1}\times\R$ with finite total mass, and that 
$\|R\mu\|_M \le \|\mu\|_M$, hence $R$ is a bounded operator from $M(\R^d)$ into  
$M(S^{d-1}\times\R)$. 
It follows from (\ref{adjoint}) that the definition (\ref{Rmu1}) coincides with (\ref{Rf}), if $\mu\in L^1(\R^d)$. 
The operator $R^*$ maps to zero all functions $\psi$ that are odd functions of $(\w,p)$; therefore we could equally well consider $R\mu$ as a linear form on the subspace of even functions in  $C_0(S^{d-1}\times\R)$. It follows that the measures $R\mu$ are all even. (A measure $\nu$, considered as a set function, is called even, if $\nu(E) = \nu(-E)$ for every Borel set $E$; in terms of the linear form this is equivalent to 
$\sca{\nu}{\vf} = \sca{\nu}{\check{\vf}}$, where $\check{\vf}$ is defined by 
$\check{\vf}(z) = \vf(-z)$.)

For sequences $\mu_k \in M(\R^d)$ we shall consider weak convergence defined by the space $C_0(\R^d)$ of test functions, 
\begin{equation}  \label{weakC0}
\lim_{k\rar\i} \sca{\mu_k}{\vf} =  \sca{\mu}{\vf}  \quad \mathrm{for\  all \ \ }   
        \vf \in  C_0(\R^d) .
\end{equation} 
In mathematical literature this is often called weak* convergence, since $M(\R^d)$ is the dual of the Banach space $C_0(\R^d)$. We will sometimes consider the analogous convergence concept with test functions in the space $C_b(\R^d)$, 
\begin{equation}     \label{weakCb}
\lim_{k\rar\i} \sca{\mu_k}{\vf} =  \sca{\mu}{\vf}  \quad \mathrm{for\  all \ \ }   
        \vf \in  C_b(\R^d) .
\end{equation} 
To distinguish those concepts we shall talk about  $C_0$-weak convergence and $C_b$-weak convergence, respectively. Occasionally we shall also consider $\mathcal D$-weak convergence; what this means should be obvious. It is obvious that  
$C_b$-weak convergence implies $C_0$-weak convergence, which in turn implies 
$\mathcal D$-weak convergence, and it is easy to see that none of those implications can be reversed. 
 Finally, note that (2.9) holds if and only if the corresponding 
  Borel measures $\mu_k$ and $\mu$ satisfy
  $\lim_{k\to\infty}\mu_k(B)=\mu(B)$ for all Borel sets $B\subset \R^d$
  for which $|\mu|(\partial B)=0$, where $\partial B$ denotes the  
  boundary of the set $B$. For probability measures this equivalence 
  is part of the well known Portmanteau theorem (see e.g.~\cite{B68}).

The definition \eqref{Rmu1} shows that the Radon transform is $C_0$-weakly continuous in the sense that $\mu_k\rar\mu$ $C_0$-weakly  implies $R\mu_k\rar R\mu$ $C_0$-weakly. The same is true if $C_0$-weakly is replaced by $C_b$-weakly.
If $\|\mu_k\|_M \le C$, then it is also true that $R\mu_k$ converges $C_0$-weakly implies  
$\mu_k$ converges $C_0$-weakly; this is essentially the content of Theorem 1$'$, see Remark 2 after Theorem~1$'$. 
If $d$ is odd a slightly stronger statement is very easy to prove as follows. 
Assume $\sca{R\mu_k}{\vf} \rar 0$ for all $\vf\in \mathcal D(\R^d)$. The formula $\psi = c R^* \p_p^{(d-1)/2} R \psi$ (see \cite{H80}), valid if $d$ is odd, helps us to write an arbitrary function $\psi \in \mathcal D(\R^d)$ in the form $\psi = R^*\vf$ with $\vf = c\, \p_p^{(d-1)/2} R\psi \in \mathcal D(S^{d-1}\times\R)$, hence \eqref{Rmu1} shows that 
$\mu_k $ tends to zero $\mathcal D$-weakly. But this implies $\mu_k \rar 0$ $C_0$-weakly, since $\mathcal D$ is dense in $C_0$ and $\|\mu_k\|_M \le C$.

For an arbitrary $L^1$-function (or measure) $g(\w,p)$ on $S^{d-1}\times\R$ it is of course not possible to define a function (measure) $p\mapsto g(\w,p)$ on $\R$ for every $\w\in S^{d-1}$. However, for any $\mu \in M(\R^d)$ the measure $R\mu\in M(S^{d-1}\times\R)$ has the special property that a measure 
$R\mu(\w,\cdot) \in M(\R)$ is well defined for every $\w\in S^{d-1}$. This is very easy to see if $\mu$ is viewed as a set function. Indeed, 
$R\mu(\w,\cdot)$ is nothing but the push-forward $\pi_{\w,*} \mu$, where $\pi_{\w}$ is the projection $\R^d\ni x \mapsto x\cdot\w \in \R$; 
here the measure $\pi_{\w,*} \mu$ is defined by $\pi_{\w,*} \mu(E) = \mu(\pi_{\w}^{-1} (E))$ for every Borel set $E\ss\R$. Similarly, looking at $\mu$ as a linear form we define $\pi_{\w,*} \mu$ as follows. First define the pullback  $\pi_{\w}^*$ on test functions by 
$\pi_{\w}^*\vf = \vf\circ \pi_{\w} \in C_b(\R^d)$ for $\vf\in C_0(\R)$. Then define the push-forward
$\pi_{\w,*} \mu$ by 
$$
\sca{\pi_{\w,*} \mu}{\vf} = \sca{\mu}{\pi_{\w}^*\vf} = \sca{\mu}{\vf(x\cdot\w)}, \quad 
    \vf\in C_0(\R) . 
$$
Thus the restriction $R\mu(\w,\cdot)$ to a particular $\w$ of the Radon transform $R\mu$ is the same as the push-forward  $\pi_{\w,*} \mu$, 
\begin{equation}      \label{Rmuw}
\sca{R\mu(\w,\cdot)}{\vf} =  \sca{\mu}{\vf(x\cdot\w)}, \quad \vf\in C_0(\R) . 
\end{equation}

From the expression (\ref{Rmuw}) we also see that $\sca{R\mu(\w,\cdot)}{\vf}$ is a continuous function of $\w$ for every $\vf\in C_0(\R)$. 

Having extended the linear form $\mu$ to $C_b(\R^d)$ as explained above we can define the 
Fourier transform $\hat{\mu}$ of $\mu\in M(\R^d)$ by 
\begin{equation}   \label{muhat}
\hat{\mu}(\xi) =  \sca{\mu}{x \mapsto  e^{-ix\cdot\xi}} . 
\end{equation}
If $\mu\in M(\R^d)$, then $\hat{\mu}$ is a uniformly continuous bounded function.

\smallskip
\noindent
{\bf Lemma 2.}
The one-dimensional Fourier transform of $R\mu(\w,p)$ with respect to $p$ for fixed $\w$, denoted $\hat{R\mu}(\w,\s)$, is related to the $d$-dimensional Fourier transform of $\mu$ by 
\begin{equation}      \label{F-slice}
\hat{R\mu}(\w,\s) = \hat \mu(\s\w)  ,  \quad \s\in\R, \ \w\in S^{d-1} .  
\end{equation}

\noindent
{\it Proof.} 
For functions in $L^1(\R^d)$ the proof consists just in interpreting an iterated integral
$\int \ldots ds\,dp$ as a multiple integral over $\R^d$. For the general case one can argue as follows. Using (\ref{Rmuw}) we see that 
$$
\hat{R\mu}(\w,\s) = \sca{R\mu(\w,\cdot)}{p\mapsto e^{-i p\s}}
    = \sca{\mu}{x\mapsto e^{-i(x\cdot\w)\s}} = \hat\mu(\s\w) ,
$$
which proves the claim. 

\smallskip

If $\mu\in M(\R^d)$ and $h\in C_b(\R^d)$ then the convolution $\mu * h$ can be defined as
$$
\mu * h (x) = \sca{\mu}{h(x - \cdot)}, \quad  x\in \R^d , 
$$
which is easily seen to be a function in $C_b(\R^d)$. 
If $h\in C_0(\R^d)$, then $\mu*h\in C_0(\R^d)$. 
If $\vf$ is a function in $\mathcal D(\R^d)$ with integral equal to $1$, then $\vf_{\e}(x) = \e^{-d}\vf(x/\e)$ tends $C_b$-weakly to the Dirac measure at the origin as $\e\rar 0$, and similarly, the family of smooth functions $\mu* \vf_{\e}$ tends $C_b$-weakly to $\mu$ as $\e\rar 0$.

It is an elementary fact  that if $\nu\in M(\R)$ then there exists a function $F(t)$ with bounded variation, defined up to a constant, such that 
\begin{equation}   \label{F(t)}
\sca{\nu}{\vf} = \int_{\R} \vf(t) dF(t) = - \int_{\R}  F(t) \vf'(t) dt  , \quad \vf\in\mathcal D(\R) . 
\end{equation}
If $F(t)$ is normalized by the requirement that $\lim_{t\rar-\i} F(t) = 0$, then 
$$
F(t) = \nu(\{s\in\R;\, s < t\})    \quad \text{for a. e.} \ \ t\in\R . 
$$
The equation (\ref{F(t)}) shows that $\nu$ is the derivative of $F$ in the distribution sense, and hence $F$ can be found from $\nu$ as a primitive function of $\nu$. 
If $\nu = R\mu(\w,\cdot)$, then  $F(p) = \mu(H_{\w,p})$, hence 
\begin{equation}      \label{dpmuh}
\sca{R\mu(\w,\cdot)}{\vf} = - \int_{\R} \mu(H_{\w,p}) \vf'(p) dp , 
 \quad  \vf \in \mathcal D(\R) ,  
\end{equation}    
which shows that the derivative in the distribution sense of the function (of bounded variation)
$p\mapsto \mu(H_{\w,p})$ is equal to the Radon transform $R\mu(\w,\cdot)$.

\section{Injectivity theorems for the Radon transform}
We now turn to the injectivity  theorems for the Radon transform on the space of measures. In the literature on the Radon transform statements like Theorems B - D  are often called support theorems. 

\smallskip
\noindent
{\bf Theorem A.}
The Radon transform is injective on  $M(\R^d)$. 

\noindent
{\it Proof.}
If $R\mu(\w,p) = 0$, then $\hat{R\mu}(\w,\s) = 0$, and by the formula $\hat{R\mu}(\w,\s) = \hat \mu(\s\w)$ it follows that $\hat\mu = 0$, hence  $\mu= 0$.

\smallskip

Let us say that a continuous function $f$ on $\R^d$ is {\it rapidly decaying at infinity} if 
\begin{equation}   \label{decay1}
f(x) = \mathcal O(|x|^{-m})  \quad  \text{as}  \quad 
   |x|\rar\i    \text{\ for every\ } m>0.
\end{equation}

To define this property for measures we choose, for arbitrary $r>1$,  a continuous function $\chi_r(x)$ on $\R^d$ such that $0\le \chi_r \le 1$, $\chi_r(x)=0$ for $|x|<r-1$, and $\chi_r(x)=1$ for $|x|>r$.
The product $\phi\mu$ of a measure $\mu\in M(\R^d)$ and a bounded continuous function $\phi$ is defined by $\sca{\phi\mu}{\vf} = \sca{\mu}{\phi\vf}$ for every $\vf\in C_0(\R^d)$.

\smallskip
\noindent
{\bf Definition 1.}
We shall say that the measure $\mu$ is  {\it rapidly decaying at infinity}, if 
\begin{equation*}
\|\chi_r \mu\|_M =\mathcal O(r^{-m})  \quad  \text{as}  \quad 
      r\rar\i \text{\ for every\ }  m>0. 
\end{equation*}

If the measure $\mu$ is defined by a continuous density $f(x)$ and $\mu$ is rapidly decaying, it is not certain that $f$ is rapidly decaying in the sense of (\ref{decay1}); in fact then $f$ does not even have to be bounded. On the other hand, any convolution of $\mu$ with a compactly supported test function must satisfy (\ref{decay1}):

\smallskip
\noindent
{\bf Lemma 3.}
If $\mu\in M(\R^d)$ is rapidly decaying at infinity and $\phi\in\mathcal D(\R^d)$, then the smooth function $\mu*\phi$ is rapidly decaying in the sense of (\ref{decay1}), that is,
\begin{equation}      \label{decay2}
|\mu*\phi(x)| = \mathcal O(|x|^{-m}) \quad  \textrm{as} \ |x|\rar \i 
    \text{\ for every\ }  m > 0 . 
\end{equation}
Moreover, every derivative of $\mu*\phi$ is rapidly decaying in the same sense. 

\noindent
{\it Proof.} 
Assume that $\phi$ is supported in the ball $\{x;\, |x|\le A\}$. If $|x| > r + A$, then $\chi_r$ is equal to $1$ on the support of $y\mapsto \phi(x-y)$, hence 
\begin{align*}
  |\mu*\phi(x)| &= |\sca{\mu}{\phi(x-\cdot)}| 
  = |\sca{\mu}{\chi_r(\cdot)\phi(x-\cdot)}| = 
  |\sca{\chi_r\mu}{\phi(x-\cdot)}|\\
  &\le \| \chi_r\mu\|_M \sup|\phi| ,
\end{align*}
which proves the first claim. Using the formula $\p^{\b}(\mu*\phi) = \mu*\p^{\b}\phi$, where $\p^{\b}$ is an arbitrary mixed derivative, we obtain the second statement.

\smallskip

If $\Om$ is an open subset of $\R^d$ we shall 
denote by $C_c(\Om)$ the set of continuous functions with compact support in $\Om$. Moreover we shall denote by $\Mloc(\Om)$ the set of linear forms on  $C_c(\Om)$ that are continuous with respect to the topology of uniform convergence on compact subsets of $\Om$, that is,
\begin{equation}\label{C_K}
  |\sca{\mu}{\vf}| \le C_K \|\vf\|  
  \quad \text{for} \ \ \vf\in C_0(\Om) \ \text{with}
  \ \supp \vf \ss K
\end{equation}	
with a constant $C_K$ depending on the compact set $K\ss\Om$. 
If (\ref{C_K}) holds with a constant $C$ independent of $K$, then the 
total mass is $\le C$ and we write $\mu\in M(\Om)$. It is clear that 
the restriction of any $\mu\in\Mloc(\Om)$ to $\Om_1\ss\Om$ with closure 
$\bar{\Om_1}\ss\Om$ must belong to $M(\Om_1)$. 
The family of non-negative measures $\mu\in \Mloc(\Om)$ is the family of
\emph{Radon measures} on $\Om$; see Chapter 7 in \cite{F99}. 
Recall that a Radon
measure on an open subset $\Om$ of $\R^d$ is a non-negative Borel measure 
$\mu$ on $\Om$ with $\mu(K)<\infty$ for every compact $K \subset \Om$.

\smallskip
\noindent
{\bf Theorem B.} 
Let $K$ be a compact, convex subset of $\R^d$, let $f$ be a continuous function on  $\R^d\sm K$ decaying at infinity faster than any negative power of $|x|$, and assume that the Radon transform $Rf(L) = 0$ for all hyperplanes disjoint from $K$. Then $f=0$ on  $\R^d\sm K$. 
More generally, let $\mu\in \Mloc(\R^d\sm K)$ be a measure that is rapidly decaying at infinity and assume that the Radon transform $R\mu$ vanishes on the open set of hyperplanes not intersecting $K$.  Then $\mu = 0$ on $\R^d\sm K$. 

\noindent
{\it Proof.}
This theorem was first proved by Helgason, see \cite{H80}. Here we will give Strichartz' short proof \cite{S82}. Approximating $\mu$ by smooth functions $f_{\e}=\mu*\vf_{\e}$, where $\vf_{\e}(x)=\e^{-d}\vf(x/\e)\in\mathcal D(\R^d)$ and $\int\vf\, dx =1$, 
and using Lemma 3 we see that the second statement follows from the first and that it is sufficient to prove the first statement for smooth functions
($f_{\e}$ is defined in the complement of the closed $\e$-neighborhood $K_{\e}=K+ \{x;\, |x|\le \e\}$ of $K$, and $Rf_{\e}(L)=0$ for all $L$ not intersecting $K_{\e}$). 
To simplify notation we give the proof first for the case $d=2$. Denote the coordinates in the plane by $(x,y)$. Fix an arbitrary line $L$ in $\R^2\sm K$ and choose coordinates such that $L$ is the $x$-axis and $K$ is contained in the halfplane $y<0$. The assumption implies that the function 
$$
G(a,b) = \int_{\R} f(x, ax + b) dx 
$$
is equal to zero for all $b\ge 0$ and all $a$ sufficiently close to $0$. Differentiating $k$ times with respect to $a$ and putting $a=0$ gives  
$$
\p_a^k G(0,b) = \int_{\R} x^k  \p_y^k f(x, b) dx 
    = \big(\frac {\p}{\p b}\big)^k \int_{\R} x^k f(x,  b) dx = 0 , \quad b \ge 0 . 
$$
The decay assumption implies that those integrals converge. 
This shows that the expression $\int_{\R} x^k f(x,  b) dx$ must be a polynomial function of  degree $k-1$ in $b$ for $b\ge 0$. But the assumption implies that this function must tend to zero as $b\rar\i$, so it must be identically zero and in particular 
$\int_{\R} x^k f(x,  0) dx = 0$. Since this is true for every $k$ it follows that 
$f(x,0) = 0$ for all $x$, that is,   $f=0$ along the line $L$. And since $L$ was arbitrary we have proved that $f=0$ outside $K$. If $d > 2$ we argue similarly assuming that $L$ is the plane $x_d=0$ and considering the function 
$G(a,b) = \int_{\R^{d-1}} f(x', x'\cdot a + b) dx'$, where $x=(x',x_d)\in\R^d$ and $a\in\R^{d-1}$.

\smallskip

In order to state Theorem C we have to formulate what it means that a measure is homogeneous of degree $\a$. If a function $f(x)$ on $\R^d\sm\{0\}$ is homogeneous of degree $\a\in\R$, that is, $f(\lb x) = \lb^{\a} f(x)$ for all $\lb> 0$, then its action on test functions satisfies
$$
\int f(x)\vf(x/\lb) dx = \lb^{d} \int  f(\lb x)\vf(x) dx =  \lb^{d+\a} \int  f(x)\vf(x) dx, 
\quad \lb > 0 
$$
with $\supp\vf\ss \R^d\sm\{0\}$.
Set $\vf_{\lb}(x) = \vf(x/\lb)$. 
Therefore a measure (or, more generally, a distribution) $\mu$ on $\R^d\sm\{0\}$ is said to be homogeneous of degree $\a$ if 
$$
\sca{\mu}{\vf_{\lb}} =  \lb^{d+\a} \sca{\mu}{\vf} \quad \text{for all} \ \ 
 \vf\in\mathcal D(\R^d\sm\{0\}) \ \ \text{and all} \ \  \lb > 0 . 
$$
Thus, for a measure $\mu$ with density $f \in L^1(\R^d)$ the definition means that $\mu$ is homogeneous of degree $\a$ if an only if the function $f$ is homogeneous of degree $\a$.

\noindent
{\bf Theorem C.} 
Let $K$ be a convex, compact set, $0\in K$, and let $\mu$ be a function or a measure on $\R^d\sm\{0\}$ that is homogeneous of degree $\a$, where $\a$ is a non-integral real number $< -d$. Assume, as in Theorem B, that $R\mu = 0$ in the set of  hyperplanes disjoint from $K$. Then $\mu=0$. 

As we shall see in Section 5, the assumption that $\a$ is non-integral cannot be omitted.

\smallskip
\noindent
{\it Proof of Theorem C.} 
It was proved in \cite {We} that any solution of $Rf(\w,p)=0$ in $|p|>1$ must be equal to an infinite sum of functions that are homogeneous of integral degrees $\le -d$.  A function that is homogeneous of non-integral degree can obviously be represented in this form only if it is identically zero. 

We will also present a self-contained proof of Theorem C using the methods of this paper. To begin with, we may assume that the set $K = \{0\}$; indeed, the Radon transform $R\mu(\w,p)$ must vanish for all $p\ne 0$, since it must be homogeneous with respect to $p$. 
Since $\a < -d$, the measure $\mu$ must have infinite mass near the origin (unless $\mu=0$), so we cannot take the Fourier transform of $\mu$ in the elementary sense. However, 
it is known that any distribution in $\R^d\sm \{0\}$ that is homogeneous of non-integral degree $\a$ can be uniquely continued to a homogeneous distribution on $\R^d$ (see Appendix, and \cite[Theorem 3.2.3]{H03}). Let us denote this distribution also by $\mu$. 
Any homogeneous distribution $f$ on $\R^d$ belongs to the Schwartz class $\mathcal S'(\R^d)$ and hence has a Fourier transform defined by $\sca{\hat f}{\vf} = \sca f{\hat{\vf}}$ for $\vf\in\mathcal D(\R^d)$, and  $\hat f\in \mathcal S'(\R^d)$. (See  \cite{H03} or any text book on distribution theory.) 
  We claim that $\hat{\mu}\in\Lloc(\R^d)$, and in fact that $\hat{\mu}$ is a continuous function. To see this, take a function $\chi\in\mathcal D(\R^d)$, equal to $1$ in some neighborhood of the origin, and write 
$\mu = \mu_0 + \mu_1$ where $\mu_0 = \chi\mu$. Then $\mu_0$ is a distribution with compact support, hence $\hat{\mu_0}$ is a $C^{\i}$ function, and $\mu_1\in M(\R^d)$, hence $\hat{\mu_1}$ is continuous. This proves the claim. 
The Radon transform $R\mu$ is a distribution on $S^{d-1}\times\R$ defined by $\sca{R\mu}{\vf} = \sca{\mu}{R^*\vf}$ for $\vf\in\mathcal D(S^{d-1}\times\R)$ (see \cite{H80}), and by assumption $R\mu = 0$ on the open set $\{(\w,p);\, p \ne 0\}$. This implies that the distribution $R\mu(\w,\cdot)$ is supported at the origin for every $\w$.  
Hence, by Theorem 2.3.4 in \cite{H03} this distribution is a linear combination of the Dirac measure at the origin and its derivatives, which means that its Fourier transform $\hat{R\mu}(\w,\s)$ is a polynomial in $\s$ with coefficients that depend on $\w$ (in fact are continuous functions of $\w$). By an extension of Lemma~2 to distributions  we know that 
$\hat{R\mu}(\w,\s) = \hat {\mu}(\s\w)$, hence 
$$
\hat {\mu}(\s\w) = \sum_0^N a_k(\w) \s^k  
$$
for some $N$ (in fact $N \le \max\{0,-\a - d + 1\}$). 
On the other hand, it is known  that the Fourier transform of a homogeneous distribution is homogeneous, in this case of non-integral degree $ -\a - d  > d - d = 0$ \cite[Theorem 7.1.16]{H03}. Since the expression on the right hand side cannot be homogeneous of non-integral degree unless all $a_k(\w)=0$, it follows that $\hat{\mu} = 0$, hence $\mu=0$ and the proof is complete.

\smallskip
\noindent
\noindent
{\it Remark.} 
In the proof above we took for granted that $R\mu(\w, \cdot)$ is a well defined distribution for every fixed $\w$. That this is true is in fact not quite obvious, but can be understood as follows. If $\mu = \mu_0 + \mu_1$, where $\mu_0$ and $\mu_1$ have the same meaning as in the proof, then $\mu_1\in M(\R)$, and hence $R\mu_1(\w,\cdot)$ is a well defined element of $M(\R)$ for every $\w$. Since $\mu_0$ has compact support we can define $R\mu_0(\w,\cdot)$ by 
\begin{equation*}
\sca{R\mu_0(\w,\cdot)}{\vf} = \sca{\mu}{\vf(x\cdot\w)},  \quad \vf\in\mathcal D(\R) , 
\end{equation*}
in analogy with (2.10). The function $x\mapsto \vf(x\cdot\w)$ does not have compact support, but since $\mu_0$ has compact support we can define $\sca{\mu_0}{\psi}$ for any $\psi\in C^{\i}(\R^d)$ as $\sca{\mu_0}{\chi\psi}$, where $\chi\in \mathcal D(\R^d)$ and $\chi = 1$ in a neighborhood of the support of $\mu_0$.

A subset $Q$ of $\R^d$ will be called a \emph{cone} if $x\in Q$ implies $\lb x\in Q$ for every $\lb > 0$. 

\smallskip
\noindent
{\bf Theorem D.}  
Let $Q$ be a closed cone such that $Q\sm\{0\}$ is contained in some open halfspace $\{x\in\R^d;\, x\cdot\w > 0\}$.
Let $K$ be a convex, compact set, let $\mu\in \Mloc(\R^d\sm K)$ and assume that $\mu$ has finite mass on sets bounded away from $K$. Assume moreover that 
$\supp\mu\ss Q$. Assume that $R\mu(L) = 0$ on the open set of hyperplanes disjoint from $K$. Then $\mu=0$ on $\R^d\sm K$.

\smallskip
This theorem is a special case of Corollary 3  in \cite{Bo92}; in the latter theorem the function (distribution) is only assumed to be rapidly decaying outside the cone $Q$, and the Radon transform is allowed to be weighted with a positive and real analytic weight function. 

\noindent
{\it Proof of Theorem  D.}
The idea of the proof is to make a projective transformation  that maps $Q$ to a compact set and thereby reduce the problem to that of Theorem B, in fact to the special case of Theorem B when $\mu$ is compactly supported. 

Write $x = (x',x_d)$, where $x'=(x_1,\ldots,x_{d-1}) \in \R^{d-1}$. We may choose coordinates so that $Q$  is the cone $\{x;\, x_d \ge \de|x'|\}$ for some $\de>0$.  
Consider the projective transformation 
\begin{equation*}
x \mapsto y = \frac x {1 + x_d} = \Psi(x) ,  \quad x\in\R^d . 
\end{equation*}
Since $y_d = x_d/(1+x_d)$ it is clear that $\Psi(Q)$ is contained in the compact set 
\begin{equation*}
\de|y'| \le y_d \le 1 . 
\end{equation*}
Writing $\ti f = f\circ \Psi^{-1}$ and $\ti L = \Psi(L)$ our Radon transform in $x$-space is transformed as follows 
\begin{equation}    \label{Rtilde}
Rf(L) = \int_L f(x) ds_x = \int_{\ti{L}} \ti f(y) J(\ti L,y) ds_y .
\end{equation}
Here $ds_x$ and $ds_y$ are the Euclidean surface measures on hyperplanes in $x$ and $y$-spaces, respectively, and $J(\ti L,y)ds_y$ is the push-forward $\Psi_*(ds_x)$. 
It is an important fact that, for an arbitrary projective transformation $\Psi$, the Jacobian $J(\ti L,y)$ factors into a product of a function depending only on the point $y$  and a function depending only on the hyperplane $\ti L$.

\smallskip
\noindent
{\bf Lemma 4.}  
The Jacobian $J(\ti L,y)$ defined by $\Psi_*(ds_x) = J(\ti L,y)ds_y$ is a non-vanishing smooth function on the manifold $Z$ of pairs $(y,\ti L)$ of points $y$  and hyperplanes $\ti L$ for which $y\in \ti L$. It can be factored so that
\begin{equation}   \label{factorability}
J(\ti L,y) = J_0(\ti L) J_1(y) .
\end{equation}

\smallskip
In the terminology introduced by Palamodov the identity \eqref{factorability} says that projective transformations are {\it factorable} with respect to the family of hyperplanes \cite[Section 3.1]{Pa3}. This property was used in an essential way in the study of Radon transforms in \cite{Pa1}, \cite{Pa2},  \cite{Bo89}, \cite{Bo92}. The factorability of projective transformations was implicit already in \cite{GGG}, where a projectively invariant Radon transform was defined, operating on sections of a certain vector bundle.

\smallskip
\noindent
{\it Sketch of proof of Lemma 4.} 
Let $\Psi$ be the mapping 
\begin{equation}
\R^d \ni x \mapsto (x,1)/\sqrt{1 + |x|^2}  \in S^{d} 
\end{equation}
from $\R^d$ onto the open upper half of the unit sphere $S^d$ in $\R^{d+1}$. Let $L$ be a hyperplane $x\cdot\w = p$ in $\R^d$, and let $\Psi(L)$ be the image of $L$ under $\Psi$, which is a $d-1$-dimensional halfsphere in $S^d$. Denote the Euclidean surface measures on $L$ and $\Psi(L)$ by $ds_L$ and $ds_{\Psi(L)}$, respectively. A straightforward calculation gives 
\begin{equation}
\frac{ds_L}{ds_{\Psi(L)}} = \frac{(1 + |x|^2)^{d/2}}{(1 + p^2)^{1/2}} ,
\end{equation}
which shows that $\Psi$ is factorable with respect to the family of hyperplanes 
(c.f.\ \cite[Lemma~1]{Bo89} and \cite[Corollary 7.5, III]{Pa2}).
An affine transformation $T$ from $\R^d$ onto itself is factorable in a trivial way, because the Jacobian $ds_L/ds_{T(L)}$ depends only on the hyperplane $L$. Now, an arbitrary projective transformation can be represented as $T\Psi^{-1}A \Psi$, where $A$ is a rotation of the sphere $S^d$ and $T$ is an affine transformation, and it is obvious that a product of factorable transformations is factorable. This completes the proof.

\smallskip
\noindent
{\it End of proof of Theorem D.}
Using Lemma 4 we can write (\ref{Rtilde}) as 
\begin{equation*}
Rf(L) = J_0(\ti L) \int_{\ti L} \ti f(y) J_1(y) ds_y .
\end{equation*}
Since $J_0 \ne 0$, the assumption that $Rf(L) = 0$ for all $L$ not intersecting $K$ now implies that the Radon transform of $\ti f J_1$ in the $y$-space vanishes on the set of all $\ti L$ not intersecting $\ti K = \Psi(K)$. Since $\ti f$ is compactly supported, Theorem B implies that $\ti f J_1 = 0$ outside $\ti K$, and since $J_1\ne 0$ it follows that $\ti f = 0$ outside $\ti K$, which in turn implies that $f = 0$ outside $K$.

\smallskip
\noindent
{\it Remark.} 
Theorem C holds without change for distributions of arbitrary order. 
Theorem B holds for distributions if rapid decay at infinity is defined for instance by \eqref{decay2} being satisfied for all $\phi\in\mathcal D(\R^d)$. Theorems A and D are valid for distributions that decay sufficiently fast at infinity for the Radon transform to be defined.

\section{ The Cram\'er-Wold theorems }
We are now ready to state several versions of the Cram\'er-Wold theorem. We begin with four versions of ``uniqueness type'', each an immediate consequence of one of the theorems above. After that we will state four analogous versions for sequences of measures.  Recall that we denote the halfspace $\{x\in\R^d;\, x\cdot\w < p\}$ by $H_{\w,p}$. 

\smallskip

\noindent
{\bf Theorem 1.}
A measure $\mu\in M(\R^d)$ is uniquely determined by $\mu(H_{\w,p})$ for almost all $(\w,p)\in S^{d-1}\times\R$. 
In other words, if $\mu(H_{\w,p}) = 0$ for almost all halfspaces $H_{\w,p} \ss \R^d$, then $\mu = 0$. 

\noindent
{\it Proof.}
As we saw above, for each $\w$ the distribution derivative of $p\mapsto \mu(H_{\w,p})$ is the Radon transform $R\mu(\w,\cdot)$ of $\mu$, evaluated at $\w$. The assertion now follows from Theorem A. 

We shall use the notation $\overline{E}$ to denote the closure of a 
subset $E\ss\R^d$. We say that $E$ is bounded away from the origin 
if $\overline{E} \subset \R^d \setminus \{0\}$.

\smallskip
\noindent
{\bf Theorem 2.}
Assume that $\mu\in \Mloc(\R^d\sm\{0\})$ is rapidly decaying at infinity and that 
$\mu(H_{\w,p}) = 0$ for almost all halfspaces $H_{\w,p}$ for which 
$\overline{H_{\w,p}} \ss \R^d\sm\{0\}$. Then $\mu = 0$. 

\noindent
{\it Proof.} 
The assumption implies that the Radon transform of $\mu$ vanishes in the open set $\{(\w,p);\, p\ne 0\}$. Application of Theorem B  with $K=\{0\}$ then proves that 
$\mu = 0$. 

\smallskip

\noindent
{\bf Theorem 3a.} 
Assume that $\mu\in \Mloc(\R^d\sm\{0\})$
 is homogeneous of non-integral degree  
$\a < -d$, and that $\mu(H_{\w,p}) = 0$ for almost all halfspaces $H_{\w,p}$ for which 
$\overline{H_{\w,p}} \ss \R^d\sm\{0\}$. Then $\mu = 0$. 

\noindent
{\it Proof.}
As in the proof of Theorem 2 the assumption implies that $R\mu = 0$ in the set 
$\{(\w,p);\, p\ne 0\}$. An application of Theorem C  with $K=\{0\}$ completes the proof. 

\smallskip
 There is a more subtle version of the previous theorem where the function 
 $p\mapsto \mu(H_{\w,p}) = a(\w,p)$ is assumed to be homogeneous, but the measure $\mu$ is not. In this case we need to assume that the measure $\mu$ is non-negative. 
Note that the closed halfspace $\overline{H_{\w,p}}$ is contained in $\R^d\sm\{0\}$ if and only if $p < 0$. 

\noindent
{\bf Theorem 3b.}
Assume that $a(\w,p)$ is a locally bounded function on $S^{d-1}\times\{p\in\R;\, p<0\}$ that is homogeneous of non-integral degree $-\b < 0$ with respect to $p$. Then there exists at most one non-negative measure 
$\mu\in \Mloc(\R^d\sm\{0\})$ with finite mass on sets bounded away from the origin, such that  
\begin{equation}   \label{thm3b} 
\mu(H_{\w,p}) = a(\w,p) \quad \text{for almost all halfspaces $H_{\w,p}$ for which 
$\overline{H_{\w,p}} \ss \R^d\sm\{0\}$.} 
\end{equation}
The measure $\mu$, if it exists, is homogeneous of degree $-d-\b$ and satisfies $R\mu = \p_p a$ in $S^{d-1}\times \{p\in\R;\,p< 0\}$. 

\smallskip
\noindent
{\it Proof.} 
Here is an outline of the proof. Using the formula $\hat{R\mu}(\w,\s) = \hat{\mu}(\s\w)$ we construct a homogeneous solution $\mu_0$ of the equation $R\mu_0 = \p_p a$; then $\mu_0$ must satisfy \eqref{thm3b}. Of course we cannot know if $\mu_0$ is non-negative, and moreover,  since we do not require $\mu_0$ to be homogeneous it is far from unique as a solution to \eqref{thm3b}. However,  $\mu_0$ is the only solution of \eqref{thm3b} that can possibly be non-negative. To prove this we shall  use the fact 
the measure $\nu = \mu -  \mu_0$, which solves the equation $R\nu = 0$ in $p \ne 0$, must be equal to a finite sum $\sum h_k$ of distributions that are homogeneous of integral order, each satisfying $Rh_k = 0$ in $p\ne 0$ (c.f.\  the proof of Theorem C). Since none of the $h_k$ can be non-negative and one of the $h_k$ must dominate in the expression  $\mu = \mu_0 + \sum h_k$ either for small or for large $|x|$, $\mu$ cannot be non-negative unless all $h_k$ vanish.

Set $b(\w,p) = \p_p a(\w,p)$ for $p <  0$ and extend $b(\w,p)$ as an even function of $(\w,p)$ on $S^{d-1}\times\R$. 
A homogeneous distribution $\mu_0$ on $\R^d$ satisfying $R\mu_0(\w,p) = b(\w,p)$ for $p\ne 0$
will now be constructed using the equation 
\begin{equation}    \label{mu0hat}
\hat{\mu_0}(\s\w) = \hat b(\w,\s) . 
\end{equation}
But the function $p\mapsto b(\w,p)$ is not integrable at the origin (unless it is identically zero), so the Fourier transform 
$\hat b(\w,\s)$ is not defined in the elementary sense. 
However, since $p\mapsto b(\w,p)$ is homogeneous of non-integral order $-\b-1$, this function
can be extended for all $\w$ uniquely to a distribution on $\R$ that is homogeneous of degree $-\b - 1$ (see Appendix!).
In this way we obtain a distribution on 
$S^{d-1} \times \R$ which we shall also denote by $b(\w,p)$, and by  $\hat b(\w,\s)$ we understand the $1$-dimensional Fourier transform of this distribution with respect to $p$. 
Note that $\hat{\mu_0}(\s\w)$ is well defined by (\ref{mu0hat}), since $b$ and $\hat b$ are even functions of $(\w,p)$ and $(\w,\s)$, respectively. The function $\s\mapsto \hat b(\w,\s)$ is homogeneous of degree $-1 - (-\b-1) = \b$, and this makes $\hat{\mu_0}$ homogeneous of degree $\b$ and locally bounded, hence $\hat{\mu_0}$ is an element of the space $\mathcal S'(\R^d)$ of tempered distributions. This implies that $\mu_0$ is a well defined distribution in $\R^d$, homogeneous of degree $-d-\b$.

Assume now that $\mu$ is a non-negative measure satisfying (\ref{thm3b}).
In order to be able to use the Fourier transform we must now prove that $\mu$ can be extended to a distribution on $\R^d$. 
Let $\mu^{\e}$ be the restriction of $\mu$ to $\R^d\sm B_{\e}$, where $B_{\e}$ is the closed ball with radius $\e$ centered at the origin.  
We claim that the norm of  $\mu^{\e}$ satisfies an estimate
\begin{equation}    \label{mueps1} 
\|\mu^{\e}\|_M \le C \e^{-\b} ,  \quad \e > 0 . 
\end{equation}
To prove this we observe that we can cover $\R^d\sm B_{\e}$ by $2d$ halfspaces of the form $H_{\w^j,-\e/d} \ss \R^d\sm\{0\}$ for suitable $\w^j$, and since $\mu\ge 0$ it then follows that 
\begin{equation}     \label{mueps2}
\|\mu^{\e}\|_M \le \sum_j \mu(H_{\w^j,-\e/d}) = \sum_j a(\w^j,-\e/d) 
\le C(\e/d)^{-\b} =  C_1 \e^{-\b} .   
\end{equation}
This is known to imply that $\mu$ can be extended  to a distribution on $\R^d$ of order $< \b + 1$  (see Appendix!). 
The extension is unique up to a distribution supported at the origin. Choose any of those extensions and denote it by $\ti{\mu}$.

We shall prove that $\ti{\mu} = \mu = \mu_0$ in $\R^d\sm\{0\}$. Set $\nu =\ti{\mu} - \mu_0$. It is clear that $R\nu(\w,p) = 0$ for $p\ne 0$. 
 As in the proof of Theorem C we can now use the formula $\hat{R\nu}(\w,\s) = \hat{\nu}(\s\w)$ together with the fact that $p\mapsto R\nu(\w,p)$ is supported at $p=0$ for every $\w$ to conclude that 
\begin{equation*}
\hat{\nu}(\s\w) = \sum_0^N c_k(\w) \s^k ,  \quad \w \in S^{d-1}, \  \s\in\R ,
\end{equation*}
where each $c_k(\w)$ is a continuous function, even if $k$ is even, odd if $k$ is odd. 
This shows that we can write
\begin{equation}    \label{sumhk1}
\mu = \mu_0 + \sum_0^N h_k ,  \quad \text{in} \ \R^d\sm\{0\} ,
\end{equation}
for some $N$, where $h_k$ are homogeneous distributions of degree $-k-d$, 
defined by $\hat{h_k}(\xi) = |\xi|^k c_k(\xi/|\xi|)$ for $\xi\in\R^d\sm\{0\}$, 
which are all mapped to zero by the Radon transform. 
Note that $h_k$ is even if $k$ is even, odd if $k$ is odd. 
Now we are going to use the assumption that $\mu\ge 0$ to prove that all $h_k$ must vanish identically. Assuming the contrary we can choose $r$ and $s$, $0\le r \le s \le N$, such that $h_r$ and $h_s$ are not identically zero and $h_k=0$ for $k\notin[r,s]$. 
For any $\vf\in\mathcal D(\R^d)$ set $\vf_{\lb}(x) = \vf(x/\lb)$. Using the homogeneity properties of $\mu_0$ and $h_k$ we now obtain
\begin{equation}  \label{sumhk}
\sca{\mu}{\vf_{\lb}}  
   = \lb^{-\b}\sca{\mu_0}{\vf} + \sum_{k=r}^s \lb^{-k} \sca{h_k}{\vf}  ,  \quad 
       \vf\in\mathcal D(\R^d\sm\{0\} . 
\end{equation}
Assume first that $\b < s$. 
Since $Rh_s(\w,p) = 0$ in $p\ne 0$ we can choose $\vf\in\mathcal D(\R^d\sm\{0\}$ such that $\vf \ge 0$ and $\sca{h_s}{\vf} < 0$. In fact,  a distribution $h$ for which $\sca h{\vf} \ge 0$ for all test functions $\vf \ge 0$ is known to be a non-negative measure, so if such a $\vf$ did not exist, $h_s$ would be a non-negative measure and hence could not have vanishing Radon transform unless it were identically zero. 
If $\lb$ is small, the term with $k=s$ dominates in  (\ref{sumhk}), hence if 
$\sca{h_s}{\vf} < 0$ we get a contradiction to $\mu\ge 0$. 
Similarly, if $\b>s>r$, then the term with $k=r$ dominates for large $\lb$, so 
if we choose $\vf$ with $\sca{h_r}{\vf} < 0$, we see that 
again $\mu$ cannot be $\ge 0$. This gives a contradiction unless all $h_k$ vanish, and the theorem is proved. 

By Theorem 3.2.4 in \cite{H03} a function $f$ in $\R^d\sm\{0\}$ that is homogeneous of {\it  integral} degree $-d-m$, $m\ge 0$, can be extended to a homogeneous distribution on $\R^d$
if and only if 
\begin{equation}   \label{horm}
\int_{|x| = 1} x^{\g} f(x) ds = 0  \quad \text{for all multi-indices $\g$ with $|\g|=m$} . 
\end{equation}
The same is true for measures in $\Mloc(\R^d\sm\{0\})$ (and even for distributions) if the condition \eqref{horm} is interpreted appropriately. 
Using this fact one can prove that the statement of Theorem~C is true under the assumption that the measure $\mu$ is even and homogeneous of degree 
$\a = -d-m$ where $m$ is an odd integer. Because if $\mu$ is even and $m$ is odd, then all the integrals \eqref{horm} must vanish, so $\mu$ must have an $\a$-homogeneous  extension whose Fourier transform satisfies $\hat{\mu}(\s\w)=a_m(\w)\s^m$. Since $m$ is odd this contradicts the assumption that $\mu$ is even, unless $\mu=0$. Theorem 3a can be extended similarly.

The assertion of Theorem 3b, finally, is true if $\b$ is an odd integer and $\w\mapsto a(\w,p)$ is even. To prove this note first that $\w \mapsto b(\w,p)$ must then also be even, and since $(\w,p) \mapsto b(\w,p)$ is even, it follows that $p\mapsto b(\w,p)$ is even.  We saw that $b(\w,p)$ is homogeneous of degree $-\b - 1$, and by assumption this is an even integer. The function  $b(\w,p)$ therefore satisfies the condition \eqref{horm} (note that $d=1$ here), hence can be extended to a homogeneous distribution on $\R$ for every $\w$. The proof can now be finished just as the proof of Theorem~3b above after we have proved that there can be no term $h_{\b}$ in \eqref{sumhk1}. In fact, $\mu$ is assumed to be even, hence each homogeneous part of \eqref{sumhk1} must be even, in particular 
$\mu_0 + h_{\b}$ must be even and $\mu_0$ is even by construction, hence $h_{\b}$ is even. But $R h_{\b} = 0$  in  $p\ne 0$ and $\b$ is odd, and we saw above that this implies that $h_{\b}$ is odd. Since $h_{\b}$ is both even and odd, it follows that $h_{\b} = 0$. This completes the proof. 

 There are similar extensions of Theorems C, 3a, and 3b where the measure $\mu$ is assumed to be odd and $\a + d$ ($\b$, respectively) is an even integer.

\smallskip

\noindent
{\bf Theorem 4.} 
Let $\mu$ be a measure in $\Mloc(\R^d\sm\{0\})$ such that 
$\mu\in M(\R^d\sm B_{\e})$ for all $\e > 0$, and let
$Q$ be a closed  cone such that $Q\sm\{0\}$ is contained in some open halfspace $\{x\in\R^d;\, x \cdot \w>0\}$. Assume moreover that $\supp\mu$ is  contained in $Q$ and 
that $\mu(H_{\w,p}) = 0$ for almost all halfspaces $H_{\w,p}$ for which $\overline{H_{\w,p}} \ss \R^d\sm\{0\}$. Then $\mu = 0$. 

\noindent
{\it Proof.}
As before we know that $R\mu = 0$ in  $\{(\w,p);\, p\ne 0\}$. 
An application of Theorem D  with $K=\{0\}$ completes the proof. 

 \medskip
 
 We shall now discuss four Cram\'er-Wold theorems for sequences of measures, 
analogous to the four theorems given above.

\noindent
{\bf Theorem 1$'$.}
Assume that $\mu_k\in M(\R^d)$  is a sequence of measures with uniformly bounded norms, $\|\mu_k\|_M \le C$, and that for almost all halfspaces $H_{\w,p}$
\begin{equation}      \label{a}
\lim_{k\rar\i} \mu_k(H_{\w,p}) = a(\w,p)    . 
\end{equation}    
Then there exists a unique measure $\mu\in M(\R^d)$ with $\|\mu\|_M \le C$ such that $\mu_k$ tends $C_0$-weakly to $\mu$, that is 
\begin{equation}      \label{thm1-prime}
\lim_{k\rar\i} \sca{\mu_k}{\vf} =  \sca{\mu}{\vf} , \quad \vf \in C_0(\R^d) .
\end{equation}     
The measure $\mu$ is characterized by $R\mu = \p_p a$ on $S^{d-1}\times\R$, where the derivative is understood in the distribution sense. If, in addition, $\lim_{k\rar\i}\|\mu_k\|_M = \|\mu\|_M$, then (\ref{thm1-prime}) holds for all $\vf\in C_b(\R^d)$. 
 
\noindent
{\it Proof.} 
Since $|\mu_k(H_{\w,p})| \le \|\mu_k\|_M \le C$ it follows from Lebesgue's theorem that 
\begin{equation}      \label{lebesgue}
\lim_{k\rar\i}  \int_{S^{d-1}} \int_{\R}  \mu_k(H_{\w,p}) \vf(\w,p) dp\, d\w 
    = \int_{S^{d-1}} \int_{\R}  a(\w,p)  \vf(\w,p) dp\,  d\w
\end{equation}	
for all $\vf \in \mathcal D(S^{d-1}\times\R)$. Using the fact that the distribution derivative of $p\mapsto  \mu_k(H_{\w,p})$ is equal to $R\mu_k(\w,\cdot)$ we obtain
\begin{align}      \label{Rmuk}
\begin{split}
\lim_{k\rar\i} \sca{R\mu_k}{\vf} 
   & = - \lim_{k\rar\i}  \int_{S^{d-1}} \int_{\R}   \mu_k(H_{\w,p})\p_p\vf(\w,p) dp\,  d\w  \\
   & = - \int_{S^{d-1}} \int_{\R}  a(\w,p) \p_p\vf(\w,p) dp\,  d\w = \sca{\p_pa}{\vf} 
   \end{split}
\end{align}	
for all $\vf\in \mathcal D(S^{d-1}\times\R)$. 
Since $\|R\mu_k\|_M \le \|\mu_k\|_M \le C$ and $\mathcal D(S^{d-1}\times\R)$ is dense in $C_0(S^{d-1}\times\R)$ it follows that (\ref{Rmuk}) holds for all 
$\vf\in C_0(S^{d-1}\times\R)$, that is, $R\mu_k$ tends $C_0$-weakly to $\p_p a$. 
Since $\|\mu_k\|_M$ is bounded we can find a subsequence $\mu_k'$ that is $C_0$-weakly convergent to some limit $\mu\in M(\R^d)$. As we have seen, this implies that $R\mu_k'\rar R\mu$ $C_0$-weakly, hence $R\mu = \p_p a$. By Theorem A this condition determines $\mu$ uniquely, hence any convergent subsequence must converge to $\mu$, so the original sequence must in fact converge to $\mu$.

To prove the last statement assume that $\lim_{k\rar\i}\|\mu_k\|_M = \|\mu\|_M$. Let $\e > 0$ and take a continuous function $\chi$ such that $0 \le \chi \le 1$, $\chi = 0$ on $|x|<r$, and $\chi = 1$ on $|x| > r+1$, with $r$  so large that $\|\chi \mu\|_M < \e$. 
Using the fact that $\|\nu\|_M \le \varliminf_{k\rar\i}\|\nu_k\|_M$ for any $\mathcal D$-weakly convergent sequence $\nu_k$ with limit $\nu$ we obtain
\begin{equation}   \label{chimuk}
\lim_{k\rar\i} \|(1-\chi)\mu_k\|_M \ge \|(1-\chi)\mu\|_M \ge \|\mu\|_M - \e . 
\end{equation}
By the assumption and by (\ref{chimuk}) we can choose $k_0$ so that 
\begin{equation*}
\|\mu_k\|_M < \|\mu\|_M + \e ,    \quad \mathrm{and}  \quad 
\|(1-\chi)\mu_k\|_M  > \|\mu\|_M - 2 \e 
\end{equation*}
for $k> k_0$. Since $\chi\ge 0$ and $1 - \chi\ge 0$ we have 
$\|\mu_k\|_M = \|\chi \mu_k\|_M + \|(1 - \chi) \mu_k\|_M$, 
hence if $k > k_0$, 
\begin{equation*}
 \|\chi \mu_k\|_M = \|\mu_k\|_M -  \|(1 - \chi) \mu_k\|_M  < \|\mu\|_M + \e 
    - (\|\mu\|_M - 2 \e) = 3 \e . 
\end{equation*}
For an arbitrary $\vf\in C_b(\R^d)$ we now write 
\begin{equation*}
\sca{\mu_k - \mu}{\vf} = \sca{\mu_k - \mu}{(1 - \chi)\vf} + \sca{\mu_k - \mu}{\chi\vf}   
\end{equation*}
and observe that the first term on the right hand side tends to zero since $(1 - \chi)\vf \in C_0(\R^d)$, and the second term can be estimated by $4\e \sup|\vf|$ since $\|\chi\mu_k\|_M < 3\e$ and $\|\chi\mu\|_M < \e$. 
The proof is complete.

\noindent
{\it Remark 1.}
Since $R\mu(\w,\cdot)$ is a well defined measure on $\R$ for every $\w$ it follows in fact that $R\mu(\w,\cdot)$ is equal to $\p_pa(\w,\cdot)$ for almost every $\w$. If we also assume that (\ref{a}) holds for almost every $p$ for  every  $\w$, then we can conclude that 
$$
R\mu(\w,\cdot) = \p_p a(\w,\cdot) 
$$
for every $\w\in S^{d-1}$. To prove this, instead of (\ref{lebesgue}) we use the fact that 
\begin{equation}      \label{leb-prime}
\lim_{k\rar\i}  \int_{\R}  \mu_k(H_{\w,p}) \vf'(p) dp 
    =  \int_{\R}  a(\w,p)  \vf'(p) dp 
\end{equation}      
for every $\w$ and all test functions $\vf(p)$ in $\mathcal D(\R)$. Then observe that the left hand side of 
(\ref{leb-prime}) is equal to 
$$
- \lim_{k\rar\i} \sca{R\mu_k(\w,\cdot)}{\vf} 
$$
by (\ref{dpmuh}), and that the right hand side is equal to $- \sca{\p_p a(\w,\cdot)}{\vf}$.

\smallskip
\noindent
{\it Remark 2.}
The first part of Theorem 1$'$, if phrased in terms of the Radon transform $R$, says essentially that $R^{-1}$ is continuous in the following sense: if $\|\mu_k\|_M \le C$ and $R\mu_k$ is $C_0$-weakly convergent, then $\mu_k$ is $C_0$-weakly convergent. 

\smallskip

If the measures $\mu_k$ are positive, then the assumption $\|\mu_k\|_M \le C$ in Theorem 1$'$ can be omitted: 

\smallskip
\noindent
{\bf Corollary 1.}
Assume that $\mu_k$ is a sequence of non-negative measures in $M(\R^d)$ and that $\lim_{k\rar\i} \mu_k(H_{\w,p})$ exists for almost all halfspaces $H_{\w,p}$ and is equal to $a(\w,p)$. Then the sequence $\mu_k$ converges $C_0$-weakly to a measure 
$\mu\in M(\R^d)$ satisfying $R\mu = \p_p a$. 

\noindent
{\it Proof.}
The sequence $\mu_k(H_{\w,p})$ must be bounded for almost every $H_{\w,p}$ by (\ref{a}), so if we cover $\R^d$ by halfspaces, $\R^d = H_1 \cup H_2$, then $\|\mu_k\|_M \le \mu_k(H_1) + \mu_k(H_2) \le C$ for all $k$ since $\mu_k \ge 0$. The assertion therefore follows from Theorem 1$'$. 

\noindent
{\it Remark.}
If the measures $\mu_k$ are not assumed to be positive measures, then the assumption $\|\mu_k\|_M \le C$ cannot be omitted. As an example, take a function $f \in \mathcal D(\R)$ with $\int f\,dx = 1$ and let $\mu_k$ be the measure with density $f_k(x) = k^2 f'(kx)$. Then $\lim_{k\rar\i}\mu_k(H) = 0$ for every halfaxis $H= (c,\pm\i)$ with $c \ne 0$, but the sequence $\mu_k$ is not $C_0$-weakly convergent, since $\sca{\mu_k}{\vf}$ tends to $-\vf'(0)$ for every continuously differentiable test function $\vf$  with compact support.

Since the case of positive measures is the most interesting in probability 
theory, we shall only consider sequences of positive measures in 
Theorems 2$'$ - 4$'$. 
Recall that the family of non-negative measures in $\Mloc(\R^d\sm\{0\})$
is the family of Radon measures on $\R^d\sm\{0\}$.
If $\vf$ is a function on $\R^d$ whose support is bounded away from the
origin, then we say that $\vf$ is supported away from the origin.

\noindent
{\bf Theorem 2$'$.} 
Let $\mu_k\in \Mloc(\R^d\sm\{0\})$ be a sequence of non-negative measures satisfying 
\begin{equation}    \label{thm2prime}
\lim_{k\rar\i} \mu_k(H_{\w,p}) = a(\w,p) 
\end{equation}
for almost all halfspaces $H_{\w,p}$ with $\overline{H_{\w,p}} \ss \R^d\sm\{0\}$, where $a(\w,p)$ is rapidly decaying in the sense that 
\begin{equation*}
a(\w,p) = \mathcal O(|p|^{-m})  \quad \text{as}\ \  p\rar -\i  \quad \text{for all} \ m .
\end{equation*}
Then there exists a unique $\mu\in \Mloc(\R^d\sm\{0\})$ such that 
\begin{equation}
\lim_{k\rar\i} \sca{\mu_k}{\vf} = \sca{\mu}{\vf} 
\end{equation}
for all $\vf\in C_b(\R^d)$ that are supported away from the origin. The measure $\mu$ satisfies $R\mu = \p_p a$ in $S^{d-1}\times \{p\in\R;\,p<0\}$ (note that only negative $p$ occur in (\ref{thm2prime})).

\smallskip
\noindent
{\it Proof.} 
Let $\mu^{\e}_k$ be the restriction of $\mu_k$ to $\R^d\sm B_{\e}$. Covering  $\R^d\sm B_{\e}$ by $2d$ halfspaces $H_{\w^j,-\e/d}$ as in the proof of Theorem 3b and choosing $k_0$ so large that 
\begin{equation*}
\mu_k(H_{\w^j,-\e/d}) < a(\w^j,-\e/d) + 1  
\end{equation*} 
for all $j$ and all $k > k_0$ we obtain 
\begin{equation*}
\|\mu^{\e}_k\|_M \le \sum_j \big(a(\w^j,-\e/d) + 1 \big) = C_{\e} .
\end{equation*}
Reasoning as in the proof of Theorem 1$'$ we can find a subsequence $\nu_k^{\e}$ that is $C_0$-weakly convergent to an element $\mu^{\e}$ in $M(\R^d\sm B_{\e})$, which satisfies 
\begin{equation}   \label{Rmu-eps}
R\mu^{\e}(\w,p) = \p_p a(\w,p)  \quad \text{in} \ \  p < -\e . 
\end{equation}
Since $\mu^{\e} \ge 0$ and $a(\w,p)$ is rapidly decaying it is easy to see that $\mu_k^{\e}$ is uniformly rapidly  decaying at infinity, and that the limit $\mu^{\e}$ is also rapidly decaying. 
It follows from Theorem~2 that $\mu^{\e}$ is uniquely determined by (\ref{Rmu-eps}),
hence the original sequence $\mu_k^{\e}$ must tend $C_0$-weakly to $\mu^{\e}$. 
Since ${\e}$ is arbitrary we obtain $\mu\in \Mloc(\R^d\sm\{0\})$, and since $\mu_k^{\e}$ is uniformly rapidly  decaying (this is of course more than we need: it suffices to observe that $\|\chi_r\mu_k^{\e}\|_M$ tends uniformly to zero as $r\rar\i$ with the notation of Definition 1) we can replace $C_0$-weak convergence by $C_b$-weak convergence, and the proof is complete.

\smallskip

\noindent
{\bf Theorem 3$'$.}
Assume that $\mu_k \in M_{\textrm{loc}}(\R^d \setminus \{0\})$ is a
sequence of non-negative measures
that satisfies
$\lim_{k\to\infty}\mu_k(H_{\omega,p})=a(\omega,p)$ for almost all
halfspaces $H_{\w,p}$ with $\overline{H_{\w,p}} \ss \R^d\sm\{0\}$, where
$p \mapsto a(\omega,p)$ is homogeneous of non-integral degree
$-\beta < 0$ for every $\omega$.
Then $\mu_k$ converges $C_b$-weakly to a measure
$\mu \in \Mloc(\R^d \sm \{0\})$ in the sense that
$\lim_{k\to\infty}\langle\mu_k,\vf\rangle=\langle\mu,\vf\rangle$
for all $\vf \in C_b(\R^d)$ supported away from the origin.
The measure $\mu$ is homogeneous of degree $-\beta-d$ and is uniquely
determined by the condition $R\mu=\p_p a$ in $S^{d-1}\times \{p\in\R;\,p<0\}$.

\smallskip
\noindent
{\it Proof.}
As above we denote 
the restriction of $\mu$ to $\R^d\sm B_{\e}$ by $\mu^{\e}$. As in the proof of Theorem~2$'$ we prove that the sequence of norms $\|\mu_k^{\e}\|_M$ is bounded by a constant $C_{\e}$. Therefore we can find a $C_0$-weakly convergent subsequence $\nu_k^{\e}$ with limit $\mu^{\e}$ satisfying $R\mu^{\e} = \p_p a$ in $p<-\e$. Since $\b$ is non-integral Theorem 3b tells us that $\mu^{\e}$ is uniquely determined by this condition. Since $\e$ is arbitrary we obtain $\mu$ in $\R^d\sm \{0\}$ such that 
$\lim_{k\rar\i}\sca{\mu_k}{\vf} = \sca{\mu}{\vf}$
for all $\vf\in C_0(\R^d)$ that are supported away from the origin. Finally, since $\mu_k\ge 0$ and $\lim_{p\rar-\i} a(\w,p) = 0$, it is clear that the total mass of $\mu_k$ outside the ball $\{x;\,|x|\le r\}$ tends uniformly to zero as $r\rar\i$, hence we obtain the same statement for $\vf\in C_b(\R^d)$ that are supported away from the origin. The proof is complete.

\smallskip
\noindent
{\it Remark.} The assumptions of Theorem 3$'$ imply that  $a(\w,p)$ must be $m$ times continuously differentiable in $\w$ if $m < \b$. To prove this we
use the fact that 
$
\sca{\p_p a(\w,\cdot)}{\psi} = \sca{R\mu(\w,\cdot)}{\psi} = \sca{\mu}{\psi(x\cdot\w)}
$
and differentiate $m$ times with respect to $\w$ observing that 
$x^{\g}\mu \in M(\R^d\sm B_{\e})$ if $|\g| \le m < \b$ since $\mu$ is homogeneous of degree $-\b - d$.

\smallskip
\noindent
{\bf Theorem 4$'$.}
Let $\mu_k$ be a sequence of non-negative measures in $\Mloc(\R^d\sm\{0\})$, such that  $\lim_{k\rar\i} \mu_k(H_{\w,p}) = a(\w,p)$ exists for almost all halfspaces 
$H_{\w,p}$ for which $\overline{H_{\w,p}} \ss \R^d\sm\{0\}$. 
Assume that there exists an open set $V\ss S^{d-1}$ such that $a(\w,p)=0$ for all $\w\in V$ and all $p<0$.
Then $\mu_k$ converges $C_0$-weakly to a  measure
$\mu \in \Mloc(\R^d \sm \{0\})$ in the sense that
$\lim_{k\to\infty}\langle\mu_k,\vf\rangle=\langle\mu,\vf\rangle$
for all $\vf \in C_0(\R^d)$ supported away from the origin. 
If $\lim_{p\rar-\i}a(\w,p) = 0$, then 
$\lim_{k\to\infty}\langle\mu_k,\vf\rangle=\langle\mu,\vf\rangle$
for all $\vf \in C_b(\R^d)$ supported away from the origin. 
The measure $\mu$ is uniquely
determined by the condition $R\mu=\partial_p a$ in $S^{d-1}\times \{p\in\R;\,p<0\}$.

\noindent
{\it Proof.} 
Again we denote by $\mu_k^{\e}$ the restriction of $\mu_k$ to $\R^d\sm B_{\e}$, and as before we prove that $\|\mu_k^{\e}\|_M \le C_{\e}$ for every $\e > 0$. Let $\nu^{\e}$ be the limit of a $C_0$-weakly convergent subsequence of $\mu_k^{\e}$; as before  $\nu^{\e}$ must satisfy the equation $R\nu^{\e} = \p_p a$ in $S^{d-1}\times \{p\in\R;\,p<-\e\}$. 
The assumptions on $a(\w,p)$ imply that $\nu^{\e}$ must be supported in a cone $Q$ satisfying the assumptions in Theorem 4.  
It now follows from Theorem 4  that $\nu^{\e}$ is uniquely determined in $\R^d\sm B_{\e}$  by the condition $R\nu^{\e} = \p_p a$, hence we obtain in this way a measure $\mu \in \Mloc(\R^d \sm \{0\})$  such that 
$\lim_{k\to\infty}\langle\mu_k,\vf\rangle=\langle\mu,\vf\rangle$
for all $\vf \in C_0(\R^d)$ supported away from the origin. To see that we can replace $\vf \in C_0(\R^d)$ by $\vf \in C_b(\R^d)$ here if $\lim_{p\rar-\i}a(\w,p) = 0$, we argue exactly as at the end of the proof of Theorem~2$'$. 

\smallskip

To formulate the next corollary we need the notation $f_{(t)}(x) = t^d f(tx)$, $t>0$, for the scaling of a function $f(x)$ on $\R^d$.  For a measure $\mu\in M(\R^d)$ the analogous operation can be defined by 
\begin{equation*}
\sca{\mu_{(t)}}{\vf} = \sca{\mu}{\vf(\cdot/t)} ,  \quad \vf\in C_0(\R^d) , 
\end{equation*}
or $\mu_{(t)}(E) = \mu(tE)$, if $\mu$ is considered as a set function.

\smallskip
\noindent
{\bf Corollary 2.} Let $\rho$ be a probability measure and
consider the family of measures
\begin{align}\label{eqcor21}
  \mu_t = t^{\beta}l(t)\rho_{(t)}, \quad t > 0,
\end{align}
where $\beta>0$ and $l$ is a positive and measurable function which
satisfies $\lim_{t\to\infty}l(\lambda t)/l(t)=1$ for every $\lambda>0$.
Assume that
\begin{align}\label{eqcor22}
  \lim_{t\to\infty}\mu_t(H_{\omega,-1})=b(\omega)
\end{align}
exists for all $\omega \in S^{d-1}$.
Assume moreover that 
\vspace{3pt}  \newline
(i)\phantom{i}  \quad $\b$ is non-integral, or  \newline 
(ii)  \quad $b(\w) = 0$ for all $\w$ in some open set $V\ss S^{d-1}$.
\vspace{3pt}  \newline  
Then $\mu_t$ converges $C_b$-weakly to some measure $\mu\in \Mloc(\R^d\sm\{0\})$ 
in the sense that 
$\lim_{t\to\infty}
\langle \mu_t,\varphi \rangle=\langle \mu,\varphi \rangle$
for all $\varphi \in C_b(\R^d)$ supported away from the origin, and $\mu$ is uniquely determined by
 $R\mu(\w,p) = \b|p|^{-\b-1}b(\w)$ in $S^{d-1}\times \{p\in\R;p<0\}$.

\noindent
{\it Proof.}
Noting that $H_{\w,\lb p} = \lb H_{\w,p}$ for $\lb>0$ we see that 
\begin{align*}
  \mu_t(H_{\omega,\lambda p})
  &=t^{\beta}l(t)\rho_{(t)}(H_{\omega,\lambda p})
  =t^{\beta}l(t)\rho_{(\lambda t)}(H_{\omega,p})\\
  &=\lambda^{-\beta}\frac{l(t)}{l(\lambda t)}
  (\lambda t)^{\beta}l(\lambda t)\rho_{(\lambda t)}(H_{\omega,p})
  =\lambda^{-\beta}\frac{l(t)}{l(\lambda t)}
  \mu_{\lambda t}(H_{\omega,p}) . 
\end{align*}
With $\lb=1/|p|$ we can conclude from (\ref{eqcor22}) that 
$\lim_{t\to\infty}\mu_t(H_{\omega,p})$ exists for for all $\w\in S^{d-1}$ and all $p<0$ and that $\lim_{t\to\infty}\mu_t(H_{\omega,p})= |p|^{-\b}b(\w)$. 
The assertion therefore follows immediately from Theorem~3$'$ if (i) holds, and from Theorem~4$'$ if (ii) holds.

\smallskip

A probability measure $\rho$ on $\R^d$ is said to be regularly varying
if there exist a $\beta>0$, a 
positive and measurable function $l$ satisfying
$\lim_{t\to\infty}l(\lambda t)/l(t)=1$ for every $\lambda > 0$
and a non-zero Borel measure $\mu$ on $\R^d\setminus \{0\}$ with 
finite mass on sets bounded away from the origin, such that 
\begin{align*}
  \lim_{t\to\infty}t^{\beta}l(t)P(tB)=\mu(B)
\end{align*}
for all Borel sets $B\subset\R^d$ bounded away from the origin
with $\mu(\partial B)=0$,
see e.g.~\cite{BDM02}, \cite{HL061}, \cite{HL062}, \cite{MS01} or \cite{R87}. 
Equivalently (see e.g.~Theorems 2.1 and 3.1 in \cite{HL062}) 
$\rho$ is regularly varying if
$\mu_t$ in \eqref{eqcor21} converges $C_b$-weakly to some non-zero 
$\mu\in M_{\textrm{loc}}(\R^d\setminus \{0\})$ as $t\to\infty$, 
in the sense that 
$\lim_{t\to\infty}\langle \mu_t,\varphi\rangle 
=\langle\mu,\varphi \rangle$ for all $\varphi\in C_b(\R^d)$ supported
away from the origin.


Hence, Corollary 2 is a characterization of regular variation for probability
measures on $\R^d$.
With assumption (i) this characterization has been shown in \cite{BDM02}.
It follows from the counterexample in Section~5 that the assumptions in 
Corollary 2 are sharp; if neither (i) nor (ii) are satisfied, then 
the conclusion need not hold.

We now discuss two applications of the characterization of regular 
variation given by Corollary 2.
The random variables considered are assumed to be defined on some 
common probability space $(\Omega,\mathcal{F},\Prob)$.

\smallskip
\noindent
{\it Random difference equations.}
Consider the random difference equation 
\begin{align}\label{rde}
  Y_n = M_nY_{n-1}+Q_n, \quad n \geq 1,
\end{align}
where $Y_n$ and $Q_n$ are $\R^d$-valued random variables and 
$M_n$ is a random $d\times d$ matrix with $\R$-valued entries.
It is assumed that the pairs $(M_n,Q_n)$, $n\geq 1$,
are independent and identically distributed.
Under weak conditions (see e.g.~\cite{K73}) the series 
\begin{align*}
  R=\sum_{k=1}^{\infty}M_1\dots M_{k-1}Q_k
\end{align*} 
converges $\Prob$-almost surely and
the probability distribution of 
$Y_n$ converges $C_b$-weakly to that of $R$, independently of $Y_0$.
If $M_1$ and $Q_1$ have non-negative entries 
and the weak (but technical) assumptions in Theorems 3 and 4 in 
\cite{K73} are satisfied, then there exists a $\beta > 0$ such that
for each $\omega \in S^{d-1}$
\begin{align*}
  \lim_{t\to\infty}t^{\beta}\Prob(\omega \cdot R > t)
\end{align*}
exists and is strictly positive for 
$\omega \in S^{d-1}_+=\{\w \in S^{d-1};\, \w_k\geq 0,\, k=1,\dots,d\}$.
Since $R$ has non-negative entries it follows that 
the limit is zero for $\omega \in -S^{d-1}_+$, so Corollary 2 implies that
the probability distribution $\rho$ given by 
$\rho(E) = \Prob(R\in E)$ is regularly varying 
with index $\beta$.

\smallskip
\noindent
{\it Domains of attraction for sums.}
Consider a sequence $\{X_k\}_{k\geq 1}$ of independent and identically
distributed $\R^d$-valued random variables. Let $X$ denote a generic 
element of the sequence and denote by $\rho$ its probability 
distribution.
It is well-known that if there exist positive constants $a_n$ and 
$\R^d$-valued constants $b_n$ such that 
the probability distribution $G_n$ of
\begin{align}\label{eqrescaledsum}
  a_n^{-1}(X_1+\dots+X_n) - b_n
\end{align}
converges $C_b$-weakly to some non-degenerate probability measure $G$, then 
$G$ is a stable distribution with characteristic exponent 
$\beta \in (0,2]$. In this case, $\rho$ is said to belong to the 
domain of attraction of $G$.
It is well known that the class of stable distributions coincides with 
the possible non-degenerate limit distributions of scaled sums of the 
type in \eqref{eqrescaledsum}.

By Theorem 4.2 in \cite{R62}, $\rho$ is in the domain of attraction of a
non-degenerate stable distribution $G$ with characteristic exponent
$\beta < 2$ if and only if $\rho$ is regularly varying with index 
$\beta$.
From Corollary 2 it follows that if $X$ takes values in $\R^d_+$, 
then $\rho$ is in the domain of attraction of a
non-degenerate stable distribution $G$ with characteristic exponent
$\beta < 2$ if and only if \eqref{eqcor22} holds with the same $\beta$
and some $l$.

\section{Counterexamples}
We need counterexamples to show  (1) that the rapid decay assumption in Theorem B cannot be omitted, (2) that the assumption that $\a$ is non-integral in Theorem C cannot be omitted, and (3) that the assumption that the cone $Q$ is contained in an open halfspace in Theorem D cannot be weakened very much. As we shall see, one sufficiently strong example meets all those requirements. 

The following simple example for dimension $d=2$, which takes care of (1) and (2) in dimension $2$, has been known for a long time (see e.g.~\cite{H80}). Let $f(x)$ be the analytic function $1/(x_1+ix_2)^2$ in $\R^2\sm\{0\}$. We claim that $\int_L f\,ds = 0$ for each line $L$ not containing the origin. In fact, the complex line integral $\int _L f(z)\,dz = \int _L z^{-2} dz$, where we have written $z = x_1 + i x_2$, is equal to zero, because  $-1/z$ is a primitive function of the integrand and it vanishes at infinity. And since $dz$ is equal to $ds$ multiplied by a non-zero complex constant along the path $L$, the claim follows. 

Similarly we can of course take $f(z) = 1/z^k$ for any $k\ge 2$. 

A closer analysis shows that the only properties of the function $f(z)=1/z^2$ that are needed here are that it is homogeneous of degree $-d = -2$, even, and has mean zero over circles centered at the origin. So, let $d$ be arbitrary $\ge 2$ and let $f(x)$ be a $C^{\i}$ function on $\R^d\sm\{0\}$, homogeneous of degree $-d$, even, and with mean zero over spheres centered at the orign. We claim that $\int_L f\,ds = 0$ for every hyperplane $L$ not containing the origin. Let $G(x)$ be the vector field 
$$
G(x) = f(x)(x_1,\ldots, x_d) .
$$
Using Euler's formula for homogeneous functions, 
$\sum x_j\p f/\p x_j = -d\, f(x)$, it is easy to see that $\divv G = 0$ in $\R^d\sm\{0\}$. Let $L_{\w,p}$ be an arbitrary hyperplane with $p\ne 0$, and consider the region in $\R^d$ that is bounded by the pair of hyperplanes $L_{\w,p}$, $L_{\w,-p}$, and the sphere $|x| = \de$, where $\de < |p|$. Since $f$ is even, $\int_{L_{\w,p}}  f\,ds = \int_{L_{\w,-p}} f\,ds$. Stokes' theorem now gives 
$2 \int_{L_{\w,p}} f\,ds = -  \int_{|x|=\de} f\, ds = 0$, which proves the claim.  
Moreover, for any mixed derivative $\p_x^{\b}$, $\b=(\b_1,\ldots,\b_d)$, the function $g=\p_x^{\b}f$ must have the same property. To see this, just observe that for given $\w$ any derivative $\p_x^{\b}$ can be written as a linear combination of derivatives of the form $D_{\w}^k D'^{\g}$, where $D_{\w}$ is the directional derivative in the direction $\w$ and $D'^{\g}$ is some derivative in the orthogonal subspace $\w^{\bot}$. Since $g=\p_x^{\b}f$ is homogeneous of degree $-d-|\b|$ we can in this way construct functions with $Rg(\w,p)=0$ for $p\ne 0$ satisfying $|g(x)|=\mathcal O(|x|^{-m})$ as $|x|\rar\i$ for arbitrarily large $m$.

Moreover, the function $f$ in the previous paragraph can be chosen with support in an arbitrarily small, open, symmetric cone $\G$ in $\R^d$.  Thus the conclusion of Theorem D may be violated if the cone $Q$ is allowed to contain an arbitrarily small conic neighborhood of a closed halfspace. 

If the dimension $d$ is $\ge 3$ we can even show that the statement of Theorem D is invalid if the cone $Q$ is assumed to be a halfspace. Indeed, take any non-zero function $h(x') = h(x_1,\ldots,x_{d-1})$ on $\R^{d-1}\sm\{0\}$, integrable at infinity, 
with $Rh(\w,p)=0$ for $\w\in S^{d-2}$ and $p\ne 0$, for instance  
$h(x')=\p_{x_1}(x_1x_2|x'|^{-d-1})$, and let $\mu$ be the measure $h(x')\de_0(x_d)$ for $(x',x_d)\in\R^d$, that is, 
\begin{equation*}
\sca{\mu}{\vf} = \int_{\R^{d-1}} \vf(x',0) h(x') dx',  \quad \vf\in C_0(\R^d), \ \                     \supp\vf\ss\R^d\sm\{0\} . 
\end{equation*}
Then $\supp\mu$ is contained in the halfspace $\{x_d\le 0\}$, and the Radon transform $R\mu$ vanishes on the set of hyperplanes not containing the origin, or expressed differently, $\mu(H)=0$ for every closed halfspace $H$ contained in $\R^d\sm\{0\}$. 
 
These examples obviously show that the assumptions of Theorems 2 - 4 are sharp in the corresponding  ways.

Finally we show that the assumptions (i) and (ii) in Corollary 2 cannot be omitted. A similar example has recently been given by Hult and Lindskog  \cite{HL061}. An advantage with the approach used here is that it makes the following 
 a very natural consequence of the examples given above. 

\smallskip
\noindent
{\bf Proposition.}
For an arbitrary integer $m\ge 1$ there exists a non-negative function $g \in L^1(\R^d)$ such that 
\begin{align}   \label{limit}  
& \lim_{t\rar\i} t^{m+d}  \int_{H_{\w,p}} g(tx) dx  
 \quad \text{exists for every halfspace $H_{\w,p} \ss \R^d\sm\{0\}$,}
\end{align}
but
\begin{align}
&  \left\lbrace 
 \begin{array} {l}   \label{nolimit} 
 \text{there exists $\vf\in\mathcal D(\R^d)$ with\ }   \text{$0\notin\supp\vf$  for which the   limit}   \\[4pt]  
  \lim_{t\rar\i} t^{m+d} \int_{\R^d} g(tx) \vf(x) dx   \quad 
\text{does not exist.} 
\end{array}
\right . 
\end{align}
The function $g$ can be chosen with support contained in an arbitrary open cone $\G$ satisfying $\G\cap(-\G)\ne\emptyset$.

\noindent
{\it Proof.}
Take $h\in C^{\i}(\{x\in\R^d;\,|x|>1\})$, not identically zero, homogeneous of degree $-d-m$ such that $\supp h\ss\G$ and $Rh(\w,p)=0$ for $|p|>1$, and set $h(x) = 0$ for $|x|<1$. To construct $g(x)$ we shall multiply $h(x)$ by a very
slowly oscillating radial function $q(x)$. This function will have the properties $|q(x)| \le 2$, 
\begin{equation}      \label{rlogr}
|\nabla q(x)| = o(|x|^{-1})  \quad \text{as} \ \  |x|\rar\i, 
\end{equation}
and for some infinite sequence of numbers $R_k$ tending to infinity
\begin{align}     \label{rings}
\begin{split}
& \text{\ \ each of the inequalities}  \quad q(x) > 1  \quad \text{and}  \quad q(x) < -1 \\
& \text{holds in infinitely many of the rings} \ \ R_k < |x| < 2R_k .
\end{split}
\end{align}
To construct such a function we can take  
$q(x) = 2\sin(\log\log |x|)$ for $|x| > e$. 
Set $g_0(x) = q(x)h(x)$ and $g(x) = g_0(x) + g_1(x)$, where 
$g_1(x) = C/ |x|^{m+d}$ for $|x|>1$ and 
$C$ is chosen so large that $g(x)\ge 0$. 
Since $t^{d+m}g_1(tx) = C|x|^{-m-d}$ for $t>1$ and $|x|>1$,
it will be enough to prove (\ref{limit}) and (\ref{nolimit}) for $g = g_0$. 
We first prove (\ref{nolimit}). Set $f_t(x) = t^{m+d}g_0(tx)$. 
By the homogeneity property of $h$ we have 
\begin{equation*}
f_t(x) = q(tx)h(x) ,  \quad \text{if\ \ } t|x| > 1 . 
\end{equation*}
Take a non-negative function $\vf\in C_0(\R^2)$, supported in a disk with radius $<1/2$ centered at the circle $|x|=2$, and so chosen that $h\ge 0$  in $\supp\vf$ and $\int h(x)\vf(x)dx = c > 0$. It follows immediately from property (\ref{rings}) that 
$\int f_t(x)\vf(x)dx = \int q(tx) h(x)\vf(x)dx$ must take values $> c$ and $< -c$  infinitely many times as $t\rar\i$.

To prove (\ref{limit}) for $g=g_0$ we write 
\begin{equation}    \label{qh}
\int_{x\cdot\w>p} q(tx)h(x) dx = 
    \int_p^{\i} \big( \int_{x\cdot\w=u} q(tx)h(x) ds\big) du . 
\end{equation}
Define a parametrization of the hyperplane $x\cdot\w=u$ by setting 
$x(y) = u\w + A_{\w}y$ for $y\in\R^{d-1}$, where $A_{\w}$ is an isometric linear map from $\R^{d-1}$ to the subspace $\{x\in\R^d;\,x\cdot\w = 0\}$.  
 To estimate the integral 
\begin{equation*}
\int_{x\cdot\w=u} q(tx)h(x) ds = \int_{\R^{d-1}} q(tx(y)) h(x(y)) dy 
\end{equation*}
we shall use the following simple estimate. If $v(y)\in C^1$ and $k(y)$ are defined on $\R^{d-1}$ and $\int k(y)dy = 0$, then 
\begin{equation*} 
|\int v(y)k(y) dy| = |\int (v(y) - v(0)) k(y) dy| \le 
   \sup_y|\nabla v(y)|  \int |y||k(y)| dy . 
\end{equation*}
With $v(y) = q(tx(y))$ and $k(y)= h(x(y))$ we have for arbitrary $\e>0$ by (\ref{rlogr}) for $u>1$ and sufficiently large $t$
\begin{equation*}
\sup_y |\nabla v(y)|  \le  {\e t}/{t u } 
  =  {\e}/{u} ,  \qquad    \text{and}
\end{equation*}
\begin{equation*}
\int_{\R^{d-1}} |y||k(y)| dy = \int_{\R^{d-1}} |y||h(x(y))|dy
   \le C_1 \int_{\R^{d-1}} \frac{|y| dy}{(u^2 + |y|^2)^{(d+m)/2}} = \frac {C_2}{u^m} . 
\end{equation*}
Hence 
\begin{equation*}
|\int_{x\cdot\w=u} q(tx)h(x) ds| \le  C_2\e/u^{m+1} , 
 \end{equation*}
which shows that the expression (\ref{qh}) is $< C_2 \e/ p^m$ if $t$ is large enough, and hence completes the proof.

\section{Appendix} 

\noindent
For the convenience of those of our readers who are not familiar with distribution theory we give here a very short proof of the fact that a function or measure on $\R^d\sm\{0\}$ that is homogeneous of non-integral degree $\a<-d$ can be uniquely extended to a homogeneous distribution on $\R^d$. As we have seen above an important consequence of this fact is that the Fourier transform (in the sense of the theory of distributions) of the extended distribution becomes available. 
The material in this section has been known since the 1950ies \cite{GS}, and is now described in many textbooks on distribution theory, and of course also in \cite{H03}.

For arbitrary $\g\in\R$ we define a function $x_+^{\g}$ in $\Lloc(\R\sm\{0\})$ by $x_+^{\g} = x^{\g}$ for $x>0$ and $x_+^{\g} = 0$ for $x<0$. If $\g>-1$  this function belongs to $\Lloc(\R)$. For any $\g \le -1$ we wish to extend $x_+^{\g}$ to a distribution on $\R$, that is, to define a distribution on $\R$ whose restriction to $\R\sm\{0\}$ is equal to $x_+^{\g}$. An easy way to solve this problem is to integrate $x_+^{\g}$ sufficiently many times, say $k$ times,  to obtain a continuous function $F$, and then define the extended distribution as the $k$:th order distribution derivative of $F$. Explicitly, choose $k$ such that $k + \g > 0$ and a continuous function
$F$ on $\R$ such that $F(x)=0$ for $x<0$ and
$F^{(k)}(x) = x^{\g}$ for $x>0$. Then define a distribution $\ti{x}_+^{\g}$ on $\R$ by 
\begin{equation}   \label{defx^g}
\sca{\ti{x}_+^{\g}}{\vf} = \sca{\p_x^{k} F}{\vf} = (-1)^{k} \sca F{\vf^{(k)}} 
    = (-1)^{k} \int_0^{\i} F(x) \vf^{(k)}(x) dx ,
 \quad  \vf\in\mathcal D(\R) . 
\end{equation}
Then $\ti{x}_+^{\g}$ is obviously a distribution on $\R$, and by partial integrations we verify that  $\sca{\ti{x}_+^{\g}}{\vf} = \sca{{x}_+^{\g}}{\vf}$ for all $\vf\in\mathcal D(\R\sm\{0\})$, that is, 
$\ti{x}_+^{\g} =  x_+^{\g}$ on $\R\sm\{0\}$.

For any $\g\in\R$ the function  $x_+^{\g}$ is homogeneous of degree $\g$ in $\R\sm\{0\}$. Moreover, if $\g$ is $< -1$ and non-integral, the definition \eqref{defx^g} is easily seen to produce a homogeneous distribution on $\R$. Because if $\g$ is non-integral, then $F$ will have the form $F(x)=c\  x_+^{\g+k}$, where $c $ is a constant depending on $\g$ and $k$, hence $F$ is homogeneous, and the distribution derivative of a homogeneous distribution is homogeneous, hence the claim is proved. 
Since any extension to $\R$ of 
${x}_+^{\g}$ can differ from $\ti{x}_+^{\g}$ only by a linear combination of the Dirac measure $\de_0$ at the origin and its derivatives, and those distributions are homogeneous of integral degrees, $\ti{x}_+^{\g}$ must be the unique extensions of  ${x}_+^{\g}$ that is homogeneous of degree $\g$ on $\R$.  
On the other hand, if $\g$ is a negative integer, then $F$ will contain a logarithmic factor, hence $F$ will not be homogeneous, so $x_+^{\g}$ is not a homogeneous distribution on $\R$, at least 
this argument does not prove that  $\ti{x}_+^{\g}$ is homogeneous. Using the definition of homogeneous distribution it is easy to check that in fact it isn't homogeneous as a distribution on $\R$.  It follows that no homogeneous extension of  ${x}_+^{\g}$ exists if $\g$ is a negative integer. 

More generally, let $f(x)$ be a continuous function on $\R\sm\{0\}$ that has at most polynomial growth as $|x|\rar 0$, that is, $|f(x)| \le C|x|^{-k}$ for $0  < |x|<1$ and some $k$. If $F$ is a $k$:th primitive of $f$, $F^{(k)}(x) = f(x)$, then $F$ is integrable up to the origin, so we can define an extension $\ti f\in \mathcal D'(\R)$ of $f$ as the $k$:th distribution derivative of $F$. The same procedure can be applied to an arbitrary measure $\mu\in \Mloc(\R\sm\{0\})$ for which the restriction $\mu_{\e}$ to $\e<|x|<1$ satisfies $\|\mu_{\e}\|_M \le C \e^{-m}$, because the second primitive of a measure on $\R$ is a continuous function.

Using spherical polar coordinates in $\R^d$ we will now use these simple arguments to construct extensions to $\R^d$ of homogeneous distributions defined in  $\R^d\sm\{0\}$. 
Let $f(x)$ be a locally integrable function on $\R^d\sm\{0\}$ that is homogeneous of non-integral degree $\g$, which we assume to be $<-d$. To construct an extension $\ti f \in\mathcal D'(\R^d)$ of $f$ we observe that we can write 
\begin{equation*}
f(r\w) = r^{\g} u(\w) , \quad  r > 0,  \ \w\in  S^{d-1} , 
\end{equation*}
for some function $u\in L^1(S^{d-1})$. Let $k$ be the smallest integer such that $k+ \g + d >0$, and choose a constant $c=c_{k,\g,d}$ such that  
$$
c\,\, ({\p}/{\p r})^k r^{k+\g+d-1} = r^{\g + d -1} .
$$
Then $G(r,\w) = c\, r_+^{k+\g+d-1} u(\w)$ is a locally integrable function on $S^{d-1}\times \R$ and we can define a distribution $\ti f$ of order $k$ on $\R^d$ by 
\begin{equation}      \label{ext-f}
\sca{\ti f}{\vf} = (-1)^k  \int_{S^{d-1}}\int_0^{\i} G(r,\w) \p_r^k \vf(r\w) dr d\w , \quad \vf \in \mathcal D (\R^d). 
\end{equation} 
This distribution must be equal to $f$ in $\R^d\sm\{0\}$, because if $\vf$ is supported in $\R^d\sm\{0\}$, we can make $k$ partial integrations with respect to $r$ in the inner integral and obtain 
$$
\sca{\ti f}{\vf} =   \int_{S^{d-1}} \int_0^{\i} r^{\g} u(\w) \vf(r\w) r^{d-1} dr d\w 
  = \int_{\R^d} f(x) \vf(x) dx . 
$$
It is easy to see that $\ti f$ satisfies $\sca{\ti f}{\vf(\cdot/\lb)} = \lb^{\g+d} \sca{\ti f}{\vf}$, which shows that $\ti f$
is homogeneous of degree $\g$. 

 It is easy to see that the same procedure can be applied if $f$ is a homogeneous measure or even a homogeneous distribution defined on $\R^d\sm\{0\}$. Similarly one can also show that a measure $\mu$ in $\Mloc(\R^d\sm\{0\})$, whose restriction $\mu_{\e}$ to $\{x\in\R^d;\, \e<|x|<1\}$ satisfies $\|\mu_{\e}\|_M \le C \e^{-m}$ for some $m$, can be extended to a distribution on $\R^d$.

  
\end{document}